\documentclass[10pt]{article}
\usepackage{geometry}
\usepackage{amssymb, amsmath,mathrsfs}
\usepackage{color}
\usepackage{ tikz, pgfplots, todonotes,ulem,amssymb, mathrsfs}

\usetikzlibrary{patterns, decorations.markings,arrows, calc}
\usepackage{framed}
\usepackage[colorlinks=true, pdfstartview=FitV, linkcolor=darkblue, citecolor=darkblue, urlcolor=darkblue]{hyperref}
\definecolor{shadecolor}{rgb}{0.98, 0.98, 0.9}
\definecolor{darkgreen}{rgb}{0.2, 0.5,  0}
\definecolor{darkblue}{rgb}{0.1,0.1,0.45}
\definecolor{red}{rgb}{0.9,0,0}
\def \S{\mathcal S_{ {Y} }}
\def \bea#1\eea {\begin{align} #1 \end{align}}

\def\&{\hspace{-15pt}&}
\def\nn{\nonumber}
\def \scr{\mathscr}
\def \div{ \mathrm {div}}

\def \AA{\mathbf A}
\def \BB{\mathbf B}

\def \CC{\mathbf C}
\def \HH{\mathbf H}
\def \bs {\boldsymbol}

\def \P {\mathbb P}

\def \eqref#1{(\ref{#1})}

\def \remove #1 { \sout{ {\color{red} #1}}}

\def \wt{\widetilde}

\def \wh{\widehat}
\newcommand{\C}{\mathbb{C}}
\renewcommand{\le}{\left}
\newcommand{\ri}{\right}
\newcommand{\R}{\mathbb{R}}
\newcommand{\Z}{\mathbb{Z}}

\newcommand{\1}{\mathbf{1}}

\renewcommand{\b}{\beta}
\renewcommand{\a}{\alpha}

\renewcommand{\d}{\mathrm d}
\renewcommand{\mod}{\,\mathrm{mod}\,}

\def\res{\mathop{\mathrm {res}}\limits_}

\def\be{\begin{equation}}
\def\ee{\end{equation}}

\def\bg{\begin{gathered}}
\def\eg{\end{gathered}}

\newtheorem{theorem}{Theorem}[section]

\newtheorem{example}[theorem]{Example}
\newtheorem{exercise}[theorem]{Exercise}

\newtheorem{lemma}[theorem]{Lemma}
\newtheorem{remark}[theorem]{Remark}
\newtheorem{problem}[theorem]{Riemann--Hilbert Problem}
\newtheorem{padeproblem}[theorem]{Pad\'e\ Problem}

\newtheorem{proposition}[theorem]{Proposition} 
\newtheorem{corollary}[theorem]{Corollary} 
 
\newtheorem{definition}[theorem]{Definition}
\newtheorem{assumption}[theorem]{Assumption}
\def\bet
{
\begin{shaded}
\begin{theorem}
}
\def\eet{\end{theorem} \end{shaded}}
\def\bd
{
\begin{shaded}
\begin{definition}
}
\def\ed{\end{definition} \end{shaded}}

\def\bp
{
\begin{shaded}
\begin{proposition}
}
\def\ep{
\end{proposition}\end{shaded}
}
\def\QED {\hfill $\blacksquare$\par\vskip 10pt}

\renewcommand{\theequation}{\arabic{section}.\arabic{equation}}

\makeatletter
\@addtoreset{equation}{section}
\makeatother

\begin{document}
\vspace{0.2cm}
\begin{center}
\begin{Large}
\textbf{Abelianization of Matrix Orthogonal Polynomials} 
\end{Large}
\end{center}

\begin{center}
M. Bertola$^{\dagger\ddagger\diamondsuit}$ \footnote{Marco.Bertola@\{concordia.ca, sissa.it\}}, 
\\
\bigskip
\begin{minipage}{0.7\textwidth}
\begin{small}
\begin{enumerate}
\item [${\dagger}$] {\it  Department of Mathematics and
Statistics, Concordia University\\ 1455 de Maisonneuve W., Montr\'eal, Qu\'ebec,
Canada H3G 1M8} 
\item[${\ddagger}$] {\it SISSA, International School for Advanced Studies, via Bonomea 265, Trieste, Italy }
\item[${\diamondsuit}$] {\it Centre de recherches math\'ematiques,
Universit\'e de Montr\'eal\\ C.~P.~6128, succ. centre ville, Montr\'eal,
Qu\'ebec, Canada H3C 3J7}
\end{enumerate}
\end{small}
\end{minipage}
\vspace{0.5cm}
\end{center}

\begin{abstract}

The main goal of the paper is to connect matrix polynomial biorthogonality on a contour in the plane with a suitable notion of scalar, multi-point Pad\'e\ approximation on an arbitrary Riemann surface endowed with a rational map to the Riemann sphere. To this end we introduce an appropriate notion of (scalar)  multi-point Pad\'e\ approximation on a Riemann surface and corresponding notion of biorthogonality of sections of the semi-canonical bundle (half-differentials). Several examples are offered in illustration of the new notions. 
\end{abstract}

\setcounter{tocdepth}{12}
\normalem
\tableofcontents
\section{Introduction}
The paper is mainly concerned  with the geometry of matrix  non-hermitian biorthogonal polynomials (BMOPs)  for a (possibly complex--valued)  weight matrix $W(z)$. These are 
 two sequences of polynomials $\{P_n(z)\}_{n\in \mathbb N},$ $\{P^\vee_n(z)\}_{n\in \mathbb N}$ of degrees at most $n$ such that 
\be
\int_\gamma P_n(z) W(z) P_m^\vee(z) \d z = \delta_{nm} \mathbf H_n,
\ee
for some matrices $\mathbf H_n$, where $\gamma$ is some chosen contour in $\mathbb C$ along which $W(z)$ is defined and such that all the moments are defined and finite. 

The notion of  MOPs and BMOPs is a non-commutative version of the standard notion of (bi)orthogonality of polynomials and the literature is vast and dates back to Krein \cite{Krein}. It has found several applications beyond approximation theory in several areas and is a staple of both approximation theory  and mathematical physics see for example  \cite{SinapAssche, Simon, Grunbaum1, Grunbaum2, Grunbaum3, Cantero, Castro}.

There have been a few recent publications \cite{Duits-Kuijlaars-aztec, Charlier, KuijlaarsGroot} where a connection was found between the matrix biorthogonality and ``scalar'' biorthogonality on  the ``spectral curve'' of $W(z)$, namely, the curve defined by $\det (y\1 - W(z))=0$.
In the quoted literature the matrix $W(z)$ is a rational function of $z$ so that the spectral curve is necessarily an algebraic Riemann surface. It is interesting to point out the the spectral curve is not really an invariant object associated with the class of biorthogonal polynomials in the following sense: in the scalar case multiplying the weight by a nonzero constant does not affect the resulting orthogonal polynomials. On the contrary, in the matrix case we can multiply $W$ on the two sides by arbitrary invertible constant matrices $W(z)\mapsto \wt W(z) = L W(z) R$, \ \ \ $L, R\in {\rm GL}_{n} (\C)$.
This transformation does affect the biorthogonal matrix polynomials, albeit in a simple way: $P_n \mapsto P_n L^{-1}$ and $P_n^\vee\mapsto R^{-1} P_n^\vee$.  However this same transformation does not preserve the spectral curve. 

This point notwithstanding, the idea of associating scalar biorthogonality on a covering curve is interesting and potentially useful in studying asymptotic properties. 

The general philosophy of associating a ``scalar'' set of data to a ``matrix'' is at the roots of much of the theory of integrable systems in several guises; for example (see \cite{BabelonBook} Ch. 3,5 for a recent review, and see also \cite{Adams1,Adams2, Bertola:EffectiveIMRN}) the isospectral deformations of a  rational Lax matrix $L$ can be realized as linear evolution on the Jacobian of its spectral curve of a line bundle $\scr L$, which encodes the whole matrix $L$ up to conjugation by constant invertible matrices. The same general idea is behind Hitchin's abelianization of Higgs bundles \cite{Hitchin}, namely the idea of associating to a pair of vector bundle $\scr V$ and Higgs field $\Phi\in H^0(End(\scr V)\otimes \mathcal K)$ a line bundle, $\scr L$,. over the spectral curve $\det (\eta\1-\Phi)=0$ embedded in $\mathcal K$. This should sufficiently motivate the wording of the title.
\paragraph{Summary of results.}

There are two main sets of results in the paper, the first necessary for the second. 
The first is contained in Section \ref{BOPs} and  describes the generalization to higher genus of the notion of biorthogonality. The notion that we generalize is that of {\it multi-point Pad\'e} approximation as in \cite{Njastad} and corresponding notion of biorthogonality along a contour $\gamma$ with respect to a weight function ${Y}:\gamma \to \C$.

This falls within the same general area as the recent \cite{PadeRS}, where  a suitable notion of Pad\'e\ approximation on Riemann surfaces was defined and shown how it is linked with a (scalar) orthogonality of sections in a sequence of line bundles of increasing degree. While the general idea is fruitful, we need to modify it to be applicable to the present situation. 

Here instead we consider not functions, but sections of appropriate twists of the semicanonical line bundle on a Riemann surface $\mathcal C$ of higher genus.  

Since these notions are less common, especially in the literature on approximation theory, we cursorily illustrate below how they apply in the familiar context of ordinary (scalar) orthogonal polynomials. 

\paragraph{Half-differentials.}
In genus zero, (non-hermitian)  orthogonal polynomials for a weight function ${Y}(z)$ on a contour $\gamma$ are polynomials $p_n(z)$ of degree  at most $n$ such that 
\be
\label{inortho}
\int_\gamma p_n (z) p_m(z) {Y}(z)\d z = \delta_{nm} h_n,
\ee
for some constants $h_n$. The (formalistic) point of view we take is that we have actually sections of the semicanonical bundle, $\sqrt{\mathcal K}$,  $\psi_n(z) = p_n(z) \sqrt{\d z}$ and the orthogonality  \eqref{inortho} reads instead
\be
\int_\gamma \psi_n (z) \psi_m(z) {Y}(z) = \delta _{nm} h_n.
\ee
While this is almost entirely an aesthetic modification in genus zero, it is this  the notion that we want to generalize to higher genus, where the modification  is of substance. Another effect of this point of view is in the counting of degrees; while a polynomial $p_n$ can be thought of as a meromorphic  function with a pole of order $n$ at infinity, the half-differential $p_n\sqrt{\d z}$ has a pole of order $n+1$. Thus the space of polynomials $\mathcal P_n$ of degree $n$ is identified with the space, $\scr P_n$,  of meromorphic sections of $\sqrt{\mathcal K}$ with a single pole of order at most $n+1$ at $\infty$.

The orthogonal polynomials are also the denominators of the Pad\'e\ approximation problem for the Weyl--Stiltjes function 
\be
W(z):= \int_\gamma \frac{{Y}(w)\d w}{z-w}.
\ee
In keeping with the above shift of point of view we can interpret the above as a section half-differential:
\be
\mathcal W(z):= \int_\gamma S(z,w) \psi_0(w) {Y}(w),
\ee
where we have introduced the {\it Szeg\"o} kernel $S(z,w) = \frac {\sqrt{\d z}\sqrt{\d w}}{z-w}$.
In \cite{PadeRS}  the Cauchy kernel was used instead to generalize the notion of Weyl-Stiltjes transform; that pivotal role is here played  by the Szeg\"o kernel. The definition of this kernel is classical \cite{FaySzego, Korotkin} and it involves an additional $g={\rm genus}(\mathcal C)$ complex parameters. 

In higher genus these extra $g$ complex parameters  (rather, $2g$ real ones) are encoded in the choice of a representation $\mathcal X: \pi_1(\mathcal C)\to U(1)$, or, which is the same, a flat line bundle.  
The notion of polynomials is thus  replaced by the notion of holomorphic sections of a line bundle obtained by tensoring the  flat bundle $\mathcal X$ with the bundle $\sqrt{\mathcal K}$ of half-differentials, and allowing poles at an increasing sequence $\scr D_n<\scr D_{n+1}$ of divisors of degree $n$.  The ``dual'' space is constructed similarly, but starting from the dual flat bundle $\mathcal X^\vee$ (i.e. the representation of $\pi_1$ given by the reciprocal multipliers). 

The notion of Szeg\"o kernel subordinated to the choice of $\mathcal X$ is described in Section \ref{higherBOPs}, where we also describe a notion of biorthogonality with respect to a pairing \eqref{pairing} of the form 
$$
 \le\langle \varphi , \varphi^\vee \ri\rangle:= \int_\gamma {Y} \varphi \varphi^\vee,
$$
where $\gamma$ is a contour, ${Y}:\gamma\to \C$ a (sufficiently smooth on $\gamma$) {\it weight} function and $\varphi, \varphi^\vee$ two sections of $\mathcal X\otimes \sqrt{\mathcal K}$, $\mathcal X^\vee \otimes \sqrt{\mathcal K}$ respectively.

The notion of biorthogonal section is given in Def. \ref{defbiortho} and we show how the resulting biorthogonal sections can be expressed in familiar determinantal form \eqref{detpsipsivee} and with a Heine--like formula \eqref{Heine}.

We then show how these sections can be characterized in terms of a Riemann--Hilbert problem \ref{RHPn} in a parallel way as it was recently done in \cite{PadeRS} and as it was done historically in the context of orthogonal polynomials in \cite{FIK}. The issue of existence and uniqueness of the solution is brought into focus in Prop. \ref{propRHP}. 

In Sec. \ref{Padepar} we describe the Pad\'e\ problem associated to these biorthogonal sections as sketched above. In Sec. \ref{secCDS} we describe the Christoffel--Darboux--Szeg\"o kernel and its properties, which generalizes the familiar notion of Christoffel--Darboux kernel for orthogonal polynomials.
The biorthogonality that we describe here is the higher genus version of the multi-point Pad\'e\ approximation, see for example \cite{Njastad}. In Section \ref{guide} we introduce the notions by the way of a guide example in the familiar setting of genus 0. \\[10pt]

The second part is contained in Sec. \ref{secMBOPs}. 
We add to the data a meromorphic function $Z:\mathcal C\to \C \P^1$ of degree $r$ and choose  the sequence of divisors $\scr D_n$ so that its support is contained in the pole divisor, $\Xi$,  of $Z$. 

The hinges of the construction are contained in Thm. \ref{Propideal} and Prop. \ref{propZ}: the main idea here is that any meromorphic  section of $\mathcal X\otimes \sqrt{\mathcal K}$ with poles only at the poles of $Z$ can be written as a linear combination of  any basis of $r$ sections of $H^0\big(\mathcal X\otimes \sqrt{K}(\Xi)\big)$ with {\it polynomial} coefficients in $z=Z(p)$. 
These polynomial coefficients are the entries of the Matrix Biorthogonal Polynomials, see Corollary \ref{corMBOPs} and the main Theorem \ref{thmMBOPs}. 

Section \ref{abelianization}  shows the inverse construction to that of Section \ref{secMBOPs}; namely we show how to associate to matrix orthogonal polynomials for a rational weight matrix $W(z)$ the set of data (flat bundle $\mathcal X$ etc.) described in Sec. \ref{secMBOPs}. This construction is essentially a reformulation of \cite{Charlier} and uses the theory of the spectral problem for rational Lax matrices \cite{Adams1, Adams2, BabelonBook, Bertola:EffectiveIMRN} in an essential way.\\[2pt]

Finally we conclude with an extensive list of examples in Section \ref{secExamples}: we offer two examples where the Riemann surface $\mathcal C$ is of genus zero. These two examples are outside (but not by much) of the scope of \cite{Charlier} because the weight matrix is not a rational function. They are the ``projection'' of the ordinary Laguerre and Hermite polynomials  by a rational map of degree $2$ to give a matrix orthogonality. See Sec. \ref{secLaguerre}, \ref{secHermite}.

We also consider the example of $\mathcal C$ being an elliptic curve (genus one) in Sec. \ref{genus1example}. In this case we can write explicitly all the  main objects (Szeg\"o kernel, etc.) in terms of classical functions (Jacobi theta functions, Weierstra\ss\ functions).  

We consider an example of ``finite'' orthogonality in Sec. \ref{sectorsion}. This means that the weight function ${Y}$  is meromorphic and the contour $\gamma$ is a closed contractible contour on the elliptic curve, containing all poles in the simply connected part. In this case the  (scalar) biorthogonality is easily seen to be of finite rank. This section is inspired by \cite{KuijlaarsGroot, Duits-Kuijlaars-aztec} where the matrix biorthogonal polynomials were used to express the statistical properties of the random tilings of the aztec diamond. 
The examples contained in Sec. \ref{sectorsion} and Sec. \ref{moregeneral}   provide a generalization
of the setup in loc. cit.. The precise relationship is explained in Remark \ref{compare}.

\paragraph{Notations.}
\begin{itemize}
\item[-] $\mathcal C$: a Riemann surface of genus $g$ with canonical bundle $\mathcal K$;
\item[-] $\mathbb J(\mathcal C)$ the Jacobian of $\mathcal C$;
\item[-] $\mathfrak A: \mathcal C \to \mathbb J(\mathcal C)$ the Abel map;
\item [-] $\Theta, \Theta_{[\bs \a,\bs \b]}$ the Riemann theta function (with characteristic $[\bs \a,\bs \b]$, respectively);
\item[-] if $p\in \mathcal C$ is a point or more generally for a divisor $\scr D$, then $\Theta(p)$ is a shorthand for $\Theta(\mathfrak A(p))$, or $\Theta(\scr D )$ stands for $\Theta(\mathfrak A(\scr D))$, respectively;
\item[-] $E(p,q)$: Klein's prime form [\cite{Fay}, pag 16].
\item[-] Given a line bundle $\mathcal L$ on $\mathcal C$ and a divisor $\scr D$ the vector space of holomorphic sections of of $\mathcal L(\scr D)$ is denoted by $H^0(\mathcal L(\scr D))$ and $h^0(\mathcal L(\scr D))$ denotes its dimension. Holomorphic sections of $\mathcal L(\scr D)$ are the same as meromorphic sections, $\varphi$, of $\mathcal L$ satisfying $\div (\varphi)\geq - \scr D$. 
\item[-] $\mathcal O$ the sheaf of holomorphic functions; $\mathcal O(\pm p)$ the sheaf of meromorphic functions with at most a pole at $p$ (at least a zero, for the case with $-$). 
\item[-] for $\varphi\in H^0(\scr L)$ its divisor is denoted by $\div(\varphi)$.
\end{itemize}

\section{Biorthogonality on a Riemann surface}
\label{BOPs}

\subsection{Guide example: (multi-point) Pad\'e\ approximations in genus $0$}
\label{guide}
The standard example of orthogonal polynomials is reviewed in the next paragraph. In the subsequent one we review the case of multi--point Pad\'e\ approximations. 
The two examples depend on a choice of strictly increasing  sequence $\scr D_n<\scr D_{n+1}$ of positive divisors of the indicated degree. 
For such a choice we define 
\be
\scr P_n := H^0(\sqrt{\mathcal K}(\scr D_{n+1})) = \scr P_n^\vee,\ \ \ \  \dim \scr P_n = n+1.
\ee
The Szeg\"o\ kernel is simply 
\be
S(z,w) = \frac {\sqrt {\d z}\sqrt{\d w}}{z-w}.
\ee
In genus zero any flat bundle $\mathcal X$ are trivial because $\C\P^1$ is simply connected; thus there is only one Szeg\"o kernel.
Given a contour $\gamma$ and ${Y}:\gamma \to \C$ a sufficiently smooth function (we tacitly assume the necessary conditions so that all integrals we write are convergent) we define the Weyl-Stiltjes section of $\sqrt{\mathcal K}$ as 
\be
\mathcal W(z) = \sqrt{\d z}\int_\gamma \frac {{Y}(w) \d w}{z-w}
\ee
The standard (multi point) Pad\'e\ approximation problem can be formulated as follows.
Consider a  fixed  spanning section $\varphi_0$ of $H^0 (\sqrt{\mathcal K}(\scr D_{1}))$. Note that in genus $0$ these have no zeros (e.g. $\varphi_0 = \sqrt{\d z}$ if $\scr D_1=\infty$). Let us denote by $\scr R$ the divisor of $\varphi_0$ (which is of the form $\scr R = -p_0 = -\scr D_1$ where $p_0$ is the single pole of $\varphi_0$, typically $\infty$). 

Then the Pad\'e\ problem can be formulated as  that of finding a section  $\mathfrak Q_{n-1}\in  H^0 (\sqrt{\mathcal K}(\scr D_{n+1}-\scr D_1))$ and a section $\psi_n\in H^0 (\sqrt{\mathcal K}(\scr D_{n+1}))$ such that 
\be
 {\varphi_0(p)\mathfrak Q_{n-1}(p)}
 - \psi_n(p)\mathcal W(p) = \mathcal O(-\scr D_n - \scr R).
\ee

\paragraph{Example: single-point case.} Let $\scr D_n = n\infty$. Then $\scr P_n= \scr P_n^\vee$ is identified with the space of polynomials of degree $\leq n$ by $\varphi = p(z) \sqrt{\d z}$ (note that $\sqrt{\d z}$ has a simple pole at $\infty$).
Since $\varphi_0(z) = \sqrt{\d z}$ it follows that  $\mathcal W(z)$ is essentially the usual Weyl-Stiltjes transforms of a weight function ${Y}(z)$ on $\gamma$ (times $\sqrt{\d z}$).
The Pad\'e\ problem \ref{padeprob} is the usual Pad\'e\ approximation whose denominators  notoriously give  orthogonal polynomials \cite{SzegoBook}.

\paragraph{Example: multi-point case.}
A slightly less standard case is when $\scr D_n$ consists of $n$ points (not all the same). To simplify suppose $\scr D_n  = \infty + \sum_{j=1}^{n-1} z_j $ with all points pairwise distinct. 
we can choose as (ordered) basis of $\scr P_n:= H^0(\sqrt{\mathcal K}(\scr D_{n+1}))$ the sections 
\be
\varphi_0=\sqrt{\d z} ,\ \ \ \varphi_\ell = \frac{\sqrt{\d z}}{z-z_\ell}, \  \  \ \ell\geq 1.
\ee
Noticing that $\scr R=  -\infty$ ($\varphi_0$ has only a simple pole and no zeros), the Pad\'e\ problem now reduces to the multi-point Pad\'e\ approximation \cite{Njastad} where 
\be
\psi_n(z) =\frac{ P_n(z)}{\prod_{\ell=1}^{n-1} (z-z_\ell)} \sqrt{\d z}
\ee
and $P_n$ a monic  polynomial of degree $n$. The section $\mathfrak Q_n$ is of the form 
\be
\mathfrak Q_n(z) = \frac {Q_{n-1}(z)}{\prod_{\ell=1}^{n-1} (z-z_\ell)} \sqrt{\d z}
\ee
with $Q_{n-1}$ a polynomial of degree $\leq n-1$. The Pad\'e\ problem then becomes 
\be
\frac{Q_{n-1}(z)}{P_n(z)} - \int_\gamma \frac {{Y}(w) \d w}{z-w} = \mathcal O \le(-2\infty -2\sum_{\ell=1}^{n-1} z_\ell \ri).
\ee
Then one obtains 
\be
Q_{n-1}(z) = \int_{\gamma} \frac {(P_n(z)-P_n(w)) {Y}(w)\d w}{z-w}
\ee
and $P_n$ is obtained by the requirement that 
\be
\sqrt{\d z}\int_\gamma \frac {P_n(w) {Y}(w)\d w}{(z-w)\prod_{\ell=1}^{n} (w-z_\ell)} = \mathcal O\le(-\infty-\sum_{\ell=1}^{n-1} z_\ell\ri)
\ee
namely the linear system of $n$ equations that embody the so-called {\it multiple orthogonality} 
\be
\int_\gamma \frac {P_n(w) {Y}(w)\d w}{\prod_{\ell=1}^{n} (w-z_\ell)}=0,
\  \ \ 
\int_\gamma \frac {P_n(w) {Y}(w)\d w}{(z_j-w)\prod_{\ell=1}^{n} (w-z_\ell)} = 0,\ j=1,\dots, n-1.
\ee
which precisely implies that $\psi_n$ is orthogonal to $\varphi_0, \dots, \varphi_{n-1}$. 
\subsection{Higher genus Riemann surfaces and biorthogonality}
\label{higherBOPs}
Let $\mathcal C$ be a smooth Riemann surface of genus $g$ and $Z:\mathcal C \to \C \mathbb P^1$ a meromorphic function (the ``projection'') of degree $r$. This allows us to think of $\mathcal C$ as a branched cover of $\C$ with $r$ sheets. We will denote by lowercase $z$ the {\it value} of $Z(p)$ for $p\in \mathcal C$ but otherwise try to keep the function $Z$ and the coordinate $z$ of the target space distinguished from each other.  We recall that the {\it ramification points} are the points $p$ of the curve $\mathcal C$ where $\d Z(p)=0$; viceversa, the {\it branch-points} are the same as the critical values of $Z$, i.e. the values of $Z(p)$ at the ramification points.

The points above  $z=\infty$ will be denoted $\infty^{(j)}$ with multiplicity $k_j$ so that $\sum k_j=r$ and the corresponding divisor is denoted by $\Xi$:
\be
\Xi= \sum k_j \infty^{(j)}.
\ee
We fix a square root $\sqrt{\mathcal K}$ of the canonical line bundle of $\mathcal C$    with $h^0(\sqrt{\mathcal K})=0$. Let $\mathcal X$ be a flat bundle given by a $U(1)$ representation (with the same symbol) of the fundamental group of $\mathcal C$. For an arbitrary (generic) positive divisor $\scr D$ an application of the Serre-Riemann-Roch theorem guarantees that 
$h^0( \mathcal X\otimes \sqrt{\mathcal K}(\scr D)) = \deg \scr D$: sections of this bundle can be thought of as half-differentials $\varphi$ with multiplicative monodromy $\mathcal X:\pi_1(\mathcal C)\to U(1)$ and with pole divisor not exceeding $\scr D$: $\div(\varphi)\geq -\scr D$. 
The dual line bundle $\mathcal X^\vee$ is the bundle constructed similarly to the above but with multiplier system $\mathcal X^\vee(\gamma) = \mathcal X(\gamma)^{-1}$,  $ \gamma \in \pi_1(\mathcal C)$.
\\[10pt]

\noindent {\bf General picture: biorthogonal ``polynomials''.}
Suppose now that $\scr D_n$ is an increasing sequence of positive divisors of degree $n$:
$
\scr D_n< \scr D_{n+1}.
$
Thus the corresponding sequences of vector spaces 
\be
\scr P_n:= H^0( \mathcal X \otimes \sqrt{\mathcal K}( \scr D_{n+1}))
=
\le\{\psi \text{ meromorphic section of }\mathcal X \otimes \sqrt{\mathcal K}:\ \ \ \div(\psi)\geq - \scr D_{n+1}\ri\}
\ee
and
$\scr P_n^\vee:= H^0( \mathcal X^\vee \otimes \sqrt{\mathcal K}( \scr D_{n+1}))$ 
 are in the same inclusion:
$$
\scr P_n\subseteq \scr P_{n+1},\qquad \qquad \scr P^\vee_n\subseteq \scr P^\vee_{n+1}, \ \ \ n\geq 0.
$$
\
The reason for the shift in the indexing is that we want to think of $\scr P_n$ as ``polynomials'' of degree $n$ and hence $\scr P_0$ should have dimension $1$.
\begin{example}
For $\mathcal C= \C\P^1$ and $\scr D_n= n\infty$ we have $\scr P_n = \mathcal P_n \sqrt{d z}$ where $\mathcal P_\ell$ denotes the vector space of polynomials of degree $\leq n$. 
\end{example}

We use the symbols $\scr P, \scr P^\vee$ to mean the direct sum of all the respective sequences.
Let $\wh \gamma\subset \mathcal C$ be a contour that  does not intersect any of the divisors $\scr D_n$, and ${Y}:\wh \gamma \to \C$ a smooth function (typically we will take ${Y}$ to be a meromorphic function on $\mathcal C$ restricted to $\wh \gamma$).
Then we can define a bilinear pairing:
\bea
\label{pairing}
\le\langle ,\ri\rangle:&  \scr P\times \scr P^\vee \to \C\\
& (\varphi,  \varphi^\vee ) \mapsto \ \le\langle \varphi , \varphi^\vee \ri\rangle:= \int_\gamma {Y} \varphi \varphi^\vee.
\eea
Note that the product $\varphi\varphi^\vee$ is an ordinary single--valued meromorphic differential.
Let us fix two sequences of sections $\{\varphi_n\}_{n\in \mathbb N}, \{\varphi^\vee
_n\}_{n\in \mathbb N}$   with $\varphi_n$ spanning $\scr P_{n} \mod \scr P_{n-1}$  (and respectively for the ``dual'' space).

The Gram-Schmidt orthogonalization process produces biorthogonal section $\psi_n, \psi_n^\vee$. We can write them explicitly as follows. Denote the ``bimoments'' by 
\be
\label{bimo}
\mu_{a,b} := \int_\gamma  {Y} \varphi_a \varphi_b^\vee, \ \ \ D_n:= \det\big[\mu_{a,b}\big]_{a,b=0}^{n-1}.
\ee
Then
\be
\label{detpsipsivee}
\psi_n :=
%\frac 1{D_n}  
\det \le[
\begin{array}{cccc|c}
\mu_{00} & \mu_{01} & \dots & \mu_{0, n-1} & \varphi_0\\
\mu_{10} & \mu_{11} & \dots & \mu_{1, n-1} & \varphi_1\\
\vdots&&&&\vdots  \\
\mu_{n1} & \mu_{n2} & \dots & \mu_{n, n-1} & \varphi_n
\end{array}
\ri]
\ ,\ \ \ \ \ 
\psi_n^\vee :=
%\frac 1{D_n}  
\det \le[
\begin{array}{ccccc}
\mu_{00} & \mu_{01} & \dots  & \mu_{0n}\\
\mu_{10} & \mu_{11} & \dots  & \mu_{1n}\\
\vdots&&&\vdots  \\
\mu_{n-1,0} & \mu_{n-1, 1} & \dots  & \mu_{n-1,n}\\
\hline 
\varphi_0^\vee &\varphi_1^\vee& \dots  & \varphi_n^\vee
\end{array}
\ri]
\ee
and they satisfy 
\be
\int_{\wh \gamma } {Y} \psi_n \psi_m^\vee = \delta_{nm} h_n.
\ee
A simple exercise shows that 
\be
{h_n} =  {D_{n+1}}{D_n}.
\ee
If $D_n \neq 0$ then $\psi_n$ ($\psi_n^\vee$) also span $\scr P_n\mod \scr P_{n-1}$, ($\scr P_n\mod \scr P_{n-1}$, respectively).
\begin{definition}
\label{defbiortho}
A section $\psi_n\in \scr P_n$ is a biorthogonal section if  $\langle \psi_n,\varphi_j^\vee\rangle=0$ for all $j=0,\dots, n-1$. Similarly a section $\psi_n^\vee \in \scr P_n^\vee$ is a biorthogonal section if $\langle  \varphi_j,\psi_n^\vee\rangle=0$, $j=0,\dots, n-1$.
\end{definition}
\begin{exercise}
Show that the notion of biorthogonality does not depend on the choices of the sequences 
$\{\varphi_n^\vee\}_{n\in \mathbb N}, \{\varphi_n\}_{n\in \mathbb N}$   provided that  $\varphi_n$ spans $\scr P_{n} \mod \scr P_{n-1}$. Hence the notion depends only on the choice of the sequence $\scr D_n$, $n\in \mathbb N$. This is familiar in the genus zero context and for the ``Laurent''orthogonal polynomials \cite{Bertola:LOPs, Gekhtman:Elementary}. In that context $\scr D_n$ have support at the two points  $z=0, \infty$ and the choice of sequence describes different families of orthogonal sections.
\end{exercise}

Note that biorthogonal sections are defined up to multiplicative nonzero scalar. 

\paragraph{Szeg\"o kernel.}
We need to introduce the notion of Szeg\"o\ kernel following\footnote{The two references above differ in their definition by a sign. We follow the first of them.} \cite{Korotkin, FaySzego}. A convenient way to think of it is as the {\it reproducing kernel} for the pair of spaces $\wp, \wp^\vee$ in the following sense:
\be
\varphi(p) = \res{q\in \scr D} S(p,q) \varphi(q),\  \ \ \varphi^\vee (q) = -\res{p\in \scr D}\varphi^\vee(p) S(p,q).
\ee
It is conceptually the equivalent of the Cauchy kernel for functions. More formally it is a kernel such that:
\begin{enumerate}
\item with respect to $p$ ($q$) is a section of $\mathcal X\otimes \sqrt{\mathcal K}(q)$ ($\mathcal X\otimes \sqrt{\mathcal K}(p)$, respectively), namely has a simple pole only on the diagonal $p=q$ and, under analytic continuation along closed loops $\gamma\in \pi_1(\mathcal C,p),\ \ \sigma\in \pi_1(\mathcal C,q)$ satisfies\footnote{The symbol $\pi_1(\mathcal C,p)$ denotes the group of homotopy classes of closed contours on $\mathcal C$ starting and ending at $p$.}
\be
S(p^\gamma,q) = \mathcal X(\gamma)S(p,q),  \ \ \ S(p,q^\sigma) =\mathcal X(\sigma)^{-1} S(p,q).
\ee
\item in any local coordinate $\zeta$ it is normalized as follows along the diagonal:
\be
S(p,q) = \frac {\sqrt{\d \zeta(p)}\sqrt{\d\zeta(q)}}{\zeta(p)-\zeta(q)} \Big(1 + \mathcal O
\big(\zeta(p)-\zeta(q)\big)\Big)
\label{Szegonorm}
\ee
\end{enumerate}
Formulas in terms of Riemann theta functions as well as important determinantal identities (``Fay's identities'') are presented in Appendix \ref{theta}.
The kernel defined here exists and is unique if and only if the flat bundle $\mathcal X$ does not lie  on the Theta divisor in the Picard variety of degree zero bundles. What these words mean concretely is  the following: choose a basis $\{a_j, b_j\}_{j=1}^g$ for the homology (the "$a,b$" cycles) and denote 
\be
{\rm e}^{2i\pi \alpha_j} = \mathcal X(a_j),\  \ \ \ {\rm e}^{2i\pi \beta_j} = \mathcal X(b_j).
\ee
The numbers $\alpha_j, \beta_j$ are defined modulo integers. Define the following point in the Jacobian $\mathbb J(\mathcal C)$ of the curve (we use the same symbol $\mathcal X$ with slight abuse of notation):
\be
\mathcal X:= \bs \beta - \bs \tau \bs \a\in \mathbb J(\mathcal C)
,\ \ \ \bs \a=(\a_1,\dots, \a_g)^t,\  \ \ \ \bs \b=(\b_1,\dots, \b_g)^t\in \R^g
\ee
Then we must have $\Theta(\mathcal X)\neq 0$, with $\Theta$ the Riemann theta function (see App. \ref{theta}). The fact that this is a sufficient condition is clear from the formula \eqref{Szegodef}. The necessity is more complicated but it follows from the general Riemann vanishing theorem for Theta functions \cite{Fay}.

\paragraph{Heine formula.}
A simple exercise using Andreief formula for alternants shows that we can write the orthogonal sections $\psi_n, \psi_n^\vee$ above in the multiple--integral form:
\bea
\psi_n(p) =\frac 1{n!}  \int_{\gamma^n} \det\le[\varphi_a(p_b)\ri]_{a,b=0}^n \det\le[\varphi_a^\vee(p_b)\ri]_{a,b=1}^{n} \prod_{j=1}^n {Y}(p_j)
\nn 
\\
\label{Heine}\psi_n^\vee(p) =\frac 1{n!}  \int_{\gamma^n} \det\le[\varphi_a(p_b)\ri]_{a,b=1}^n \det\le[\varphi_a^\vee(p_b)\ri]_{a,b=0}^{n} \prod_{j=1}^n {Y}(p_j)
\eea
where $p_0=p$ (note that one of the two alternants in each integral is of size $n+1$).
\paragraph{Riemann--Hilbert problem.}
For any $\varphi\in \scr P, \varphi^\vee\in \scr P^\vee $ we define 
\be
\S[\varphi](p):= -\frac 1{2i\pi}\int_\gamma S(p,q) \varphi(q) {Y}(q),\qquad
\S[\varphi^\vee](p):= -\frac 1{2i\pi}\int_\gamma  \varphi^\vee(q)S(q,p) {Y}(q).
\ee
Note that, by the Sokhotskii-Plemeljj formula,
\bea
\S[\varphi](p_+)-\S[\varphi](p_-)= \varphi(p) {Y}(p),
\qquad\S[\varphi^\vee ](p_+)-\S[\varphi^\vee](p_-)=-\varphi^\vee(p) {Y}(p),\cr \ p\in \gamma.
\eea
\begin{lemma}
\label{rhovan}
Let $\psi_n$ be the $n$-th biorthogonal section for the pairing \eqref{pairing}.  Then $\rho_n:= \S[\psi_n]$ has the property 
\be
\div(\rho_n)\geq \scr D_n.
\ee
Viceversa, if $\psi_n \in \scr P_{n}$ is such that $\rho_n=\S[\psi_n]$ vanishes at $\scr D_n$ then $\psi_n$ is a biorthogonal section. 
Similar statement holds for $ \psi_n^\vee$. 
\end{lemma}
{\bf Proof.}
Recall that $\div(\psi_n)\geq - \scr D_{n+1}$ and that $\scr  D_{n+1}> \scr D_{n}$. 
We proceed by induction on $n$. The divisor $\scr D_2 = p_1+p_2$ consists of only two points (possibly repeated)  and hence $\psi_1$ has at most two simple poles or a double pole at $\scr D_2$. Let $\varphi_0^\vee$ be a spanning section of $\scr P_0^\vee$; namely, it has a necessarily simple pole at $\scr D_1=p_1$ since by assumption $h^0(\mathcal X^\vee\otimes \sqrt{\mathcal K})=0$. 
Then, by the Cauchy theorem 
\be
\res{} \varphi_0^\vee \rho_1 =\int_\gamma \varphi_0^\vee\le( (\rho_1)_+-(\rho_1)_-\ri) = 
 \int_\gamma \varphi_0^\vee  \psi_1 {Y},\label{rhopsi}
\ee
where in the last step we have used the Sokhotskii-Plemeljj formula and $\res{}$ denotes the sum over all residues (outside $\gamma$) of the expression. In this case it is just the residue at $p_1$ where $\varphi_0^\vee$ has a pole. 
Now $\psi_1$ is a biorthogonal section if and only if \eqref{rhopsi} vanishes. But this means that $\rho_1$ must vanish at $p_1=\scr D_1$. This concludes the initial step of induction. 

\noindent{\bf Induction step.}
Similar considerations show that, for $j=0,\dots n-1$ 
\be
\res{\scr D_n} \varphi_j^\vee \rho_n = \int_\gamma \varphi_j^\vee \psi_n {Y}.
\label{rhopsin}
\ee
By induction, $\div(\rho_n)\geq \scr D_{n-1}$ so there is only one residue to consider at $\scr D_{n}-\scr D_{n-1} = p_\star$.  If this point appear with multiplicity $1$ in $\scr D_n$ then  necessarily $\varphi_{n-1}^\vee$ has a simple pole at $p_\star$ and then the vanishing of \eqref{rhopsin} implies that $\rho_n$ must vanish as well.
	
If $p_\star$ has multiplicity $k$ in  $\scr D_{n}$ (and hence multiplicity $k+1$ in $\scr D_{n+1}$)  then we can find $k+1$ sections in $\scr P_{n-1}$ with poles of orders $1,2,\dots, k+1$ at $p_\star$. The induction hypothesis implies that $\rho_n$  already vanishes of order $k$ at $p_\star$; then the pairing with a section with a pole of order $k+1$ at $p_\star$ promptly implies the vanishing of order $k+1$. \QED

Consider the two matrices 
\be
\Upsilon_n(p) := \le[\begin{array}{cc}
\psi_n & \S[\psi_n]\\
\psi_{n-1} & \S[\psi_{n-1}]
\end{array}\ri],\ \ 
\Upsilon^\vee_n(p) := \le[\begin{array}{cc}
 \S[\psi^\vee_n] & \S[\psi^\vee_{n-1}]\\
\psi^\vee_n &\psi^\vee_{n-1}
\end{array}\ri].
\ee
The Lemma \ref{rhovan} shows that their divisor properties (outside of  $\gamma$) are
\be
\Upsilon_n = \le[
\begin{array}{cc}
\mathcal O(\scr D_{n+1}) & 
\mathcal O(-\scr D_{n}) 
\\
\mathcal O(\scr D_{n}) & 
\mathcal O(-\scr D_{n-1}) 
\end{array}\ri],\ \ \ \ 
\Upsilon_n^\vee = \le[
\begin{array}{cc}
\mathcal O(-\scr D_{n}) 
 & 
\mathcal O(-\scr D_{n-1}) 
\\
\mathcal O(\scr D_{n+1}) & 
\mathcal O(\scr D_{n})
\end{array}\ri].
\label{On}
\ee
Moreover they satisfy the boundary value conditions ({\it jump} condition) 
\be
\label{jump}
\Upsilon_{n}(p_+) = \Upsilon_{n}(p_-) \le[
\begin{array}{cc}
1 & {Y}(p)\\
0 & 1
\end{array}
\ri] ,\qquad
\Upsilon_{n}^\vee(p_+) =\le[
\begin{array}{cc}
1 & -{Y}(p)\\
0 & 1
\end{array}
\ri] \Upsilon_{n}^\vee(p_-)  ,\qquad p\in \gamma,
\ee
where the subscript $p_\pm$ denotes the boundary values at $p\in \gamma$ from the left ($+$) or the right ($-$). 
Generically the $(1,1)$ entry of $\Upsilon_n$ is a section of $\scr P_n$ but not in $\scr P_{n-1}$ (i.e., it really has  an ``extra pole''). Similarly, $\S[\psi_{n-1}]$ vanishes at $\scr D_{n-1}$ but not at $\scr D_{n}$. 
These two conditions specify a Riemann--Hilbert problem on $\mathcal C$
\begin{problem}
\label{RHPn}
Find a $2\times 2$ matrix $\Upsilon_n$ ($\Upsilon_n^\vee$) of local meromorphic sections   of $\mathcal X\otimes \sqrt{\mathcal K}$ ($\mathcal X^\vee \otimes \sqrt{\mathcal K}$, respectively) satisfying the two conditions \eqref{On} and \eqref{jump}. 
Moreover we require that the $(1,1)$ entry is $\mathcal O(\scr D_{n+1})$ but not $\mathcal O(\scr D_n)$ and the $(2,2)$ entry is $\mathcal O(-\scr D_{n})$ but not $\mathcal O(-\scr D_{n+1})$
\end{problem}
Note that if $\Upsilon_n$ solves the RHP \eqref{RHPn} we are still free to multiply it on the left by a diagonal constant invertible matrix. In general we can dispose of this freedom if we choose some convenient normalization. 
If we have chosen the bases $\varphi_\ell, \varphi_\ell^\vee$ then this normalization can be fixed by requiring that 
\be
\label{normal}
(\Upsilon_n)_{11}(p) = \varphi_n(p) + \sum_{\ell=0}^{n-1} c_\ell \varphi_\ell(p),\ \qquad
\res{\scr D_{n+1}} \varphi^\vee _{n} (\Upsilon_{n})_{22} = 1.
\ee

We also observe that the product $\Upsilon_n \Upsilon_n^\vee$ is  a matrix of  differentials without discontinuity across $\gamma$, with the $(1,1),(2,2)$ entries being holomorphic differential and the $(1,2), (2,1)$ entries having poles at $\scr D_{n+1}-\scr D_{n-1}= p_n+ p_{n+1}$. 

Note also that for $n=0$ we have 
\be
\Upsilon_0 = \le[
\begin{array}{cc}
\varphi_0  & \S[\varphi_0]\\
0 & 0
\end{array}
\ri]
\ee
We now discuss uniqueness.
\bp
\label{propRHP}
The solution $\Upsilon_n$ of the RHP \ref{RHPn} (and  similarly for  $\Upsilon_n^\vee$) with the normalization \eqref{normal} exists and is unique if and only if the determinant 
\be
D_n= \det \bigg[\mu_{a,b}\bigg]_{a,b=0}^{n-1},\ \ \ \ \ \mu_{a,b} = \int_\gamma \varphi_a^\vee \varphi_b {Y}
\ee
does not vanish. In that case the $(1,1)$ entry is the biorthogonal section $\frac {\psi_n}{D_n}$ ($\frac {\psi_n^\vee}{D_n}$, respectively). 
\ep
{\bf Proof.} 
%\begin{ocg}{OCG 1}{ocg1}{0}
Denote provisionally $\psi = (\Upsilon_n)_{11}$ and $\rho = (\Upsilon_{n})_{22}$.
The jump condition \eqref{jump} implies that $\psi$ extends analytically across $\gamma$ to a meromorphic section of $\mathcal X\otimes \sqrt{\mathcal K}$. The growth conditions \eqref{On} imply that $\div(\psi)\geq - \scr D_{n+1}$ and hence $\psi\in \scr P_n$ by definition. The normalization \eqref{normal} implies that we can write 
\be
\psi (p) = \varphi_n(p) + \sum_{\ell=0}^{n-1} c_\ell \varphi_\ell(p).
\ee
Now the condition that $\S[\psi]$ (in the $(1,2)$ entry of $\Upsilon_n$) is $\mathcal O(-\mathcal D_{n+1})$ means that it vanishes at $\scr D_{n+1}$; but this is equivalent to the orthogonality of $\psi$ to all $ \varphi_{0}, \dots, \varphi_{n-1}$ by Lemma \ref{rhovan}. This condition in turn translates to the linear system;
\be
\le[
\begin{array}{ccc}
\mu_{00} & \dots & \mu_{0,n-1}\\
\vdots && \vdots\\
\mu_{n-1,0} & \dots & \mu_{n-1,n-1}
\end{array}
\ri] \le[
\begin{array}{c}
c_0\\ \vdots \\ c_{n-1}
\end{array}
\ri] 
=\le[\begin{array}{c}
\mu_{0n}\\ \vdots \\ \mu_{n-1,n}
\end{array}
\ri] 
\label{sysc}
\ee
which is then uniquely solvable if and only if the determinant $D_{n} \neq 0$.
Now consider the $(2,2)$ entry $\rho$ and let $\wt \psi $ denote the $(2,1)$ entry of $\Upsilon_n$. 
The jump condition implies that $\rho = \S[\wt \psi]$ and the growth condition implies that $\wt \psi\in \scr P_{n-1}$. 
Thus we can write
\be
\wt \psi(p) = \sum_{\ell=0}^{n-1} \beta_\ell \varphi_\ell(p). 
\ee
The vanishing condition of $\rho$ at $\scr D_{n}$ implies that $\wt \psi$ is an orthogonal section. The normalization condition reads
\be
\delta_{\ell n}= \res{\scr D_{n}} \varphi_\ell^\vee \rho = \int_{\gamma} \wt \psi \varphi^\vee_\ell {Y},\ \ \ 0\leq \ell\leq n.
\ee
This is equivalent to the following system:
\be
\le[
\begin{array}{ccc}
\mu_{00} & \dots & \mu_{0,n-1}\\
\vdots && \vdots\\
\mu_{n-1,0} & \dots & \mu_{n-1,n-1}
\end{array}
\ri] \le[
\begin{array}{c}
\beta_0\\ \vdots \\ \beta_{n-1}
\end{array}
\ri] 
=\le[\begin{array}{c}
0\\ \vdots \\ 1
\end{array}
\ri] 
\label{sysbeta}
\ee
Then $D_n\neq 0$ if and only if this system \eqref{sysbeta} is uniquely solvable. 
%\end{ocg}
\QED
As noted in \cite{PadeRS}, in genus zero the matrix of bimoments  $\mu_{ab}$ is a Hankel matrix and then the vanishing $D_n=0$ precludes the compatibility of the two systems \eqref{sysc}, \eqref{sysbeta}. In higher genus one may  have $D_n=0$ and yet still have compatibility of the two systems \eqref{sysc}, \eqref{sysbeta} in some non-generic situations. 

The RHPs \ref{RHPn} characterize (normalized) biorthogonal sections  by taking the $(1,1)$ entry of the solutions. The existence of the Riemann--Hilbert characterization of the (normalized) biorthogonal sections is a potentially useful tool to investigate asymptotic issues; this would require to develop a higher genus analog of the Deift--Zhou  steepest descent method already used successfully in the ordinary (genus zero) orthogonality \cite{DKMVZ}. 

\subsubsection{Pad\'e\ approximation}
\label{Padepar}
In a similar vein as \cite{PadeRS} we show that the biorthogonal sections $\psi_n, \psi_n^\vee$ are the denominators of a suitable Pad\'e--like pair of approximation problems for the two Weyl ``functions''
\be
\mathcal W(p) =  \int_\gamma S(p,q) \varphi_0^\vee(q) {Y}(q).\ \ \ \ 
\mathcal W^\vee(p) = \int_\gamma \varphi_0^\vee(q) S(q,p)  {Y}(q),\ \ \ \ 
\ee
They are holomorphic  sections of $\mathcal X\otimes \sqrt{\mathcal  K}$ and $\mathcal X^\vee\otimes \sqrt{\mathcal  K}$ on $\mathcal C\setminus \gamma$, with a jump discontinuity across $\gamma$.
We will consider the Pad\'e\ problem for $\mathcal W$, the one for $\mathcal W^\vee$ being formulated by simply swapping the roles of $\mathcal X, \mathcal X^\vee$. 
\begin{padeproblem}
\label{padeprob}
Given fixed  spanning sections $\varphi_0,\varphi_0^\vee$ of $H^0 (\mathcal X \otimes \sqrt{\mathcal K}(\scr D_{1}))$, $H^0 (\mathcal X^\vee \otimes \sqrt{\mathcal K}(\scr D_{1}))$, respectively, let $\scr R= \div(\varphi_0)$,  $\scr R^\vee= \div(\varphi_0^\vee)$ be their divisors (of degree $g-1$).

Find a section  $\mathfrak Q_{n-1}\in  H^0 (\mathcal X^\vee \otimes \sqrt{\mathcal K}(\scr D_{n+1}-\scr D_1)$ and a section $\psi_n\in H^0 (\mathcal X^\vee \otimes \sqrt{\mathcal K}(\scr D_{n+1})$ such that 
\be
\label{pade}
 {\varphi_0^\vee(p)\mathfrak Q_n(p)}
 - \psi_n(p)\mathcal W^\vee(p) = \mathcal O(-\scr D_n - \scr R^\vee).
\ee
\end{padeproblem}
Since the degree of $\scr D_n+\scr R^\vee$ is $n+g$, we have no choice but set 
\be
\varphi_0^\vee(p) \mathfrak Q_n(p) = \int_{\gamma} 
\bigg(\psi_n(p)\varphi_0^\vee (q) S(q,p)-\varphi_0^\vee(p) S(p,q)\psi_n(q)\bigg){Y}(q), 
\ee
and observe that the expression is a meromorphic differential without discontinuity on $\gamma$. 
Now \eqref{pade} becomes a constraint on the remainder 
\be
\varphi_0^\vee (p) \int_\gamma \psi_n(q) {Y}(q) S(p,q) = \mathcal O(-\scr D_n - \scr R^\vee)
\ee
We can then divide by $\varphi_0$ and obtain 
\be
\rho_n(p):= \int_\gamma \psi_n(q) {Y}(q) S(p,q) = \mathcal O(-\scr D_n)
\ee
As we have seen in Lemma \ref{rhovan} this is the requirement that $\psi_n$ is a biorthogonal section.
\paragraph{Orthogonality.}
The spaces $\scr P, \scr P^\vee$ are different spaces simply because $\mathcal X\neq \mathcal X^\vee$. However,  if $\mathcal X$ takes values in $\Z_2=\{1,-1\}$ then $\mathcal X^\vee=\mathcal X$ and the pairing is actually a non-hermitian inner product. In this case one promptly sees that the matrix of bimoments $\mu_{a,b}$ \eqref{bimo} is symmetric  and  $\psi_n = \psi_n^\vee$.

\subsection{Christoffel-Darboux-Szeg\"o reproducing kernels}
\label{secCDS}
\bd
The Christoffel--Darboux-Szeg\"o (CDS) kernel is the expression 
\be
\label{CDS}
K_n(p,q):= \sum_{j=0}^{n-1} \frac 1{h_j} \psi_j(p) \psi_j^\vee(q) .
\ee
\ed
It is immediate to verify that this kernel is the projection on the subspaces  $\scr P_{n-1}, \scr P^\vee_{n-1}$ in the following sense
\be
\int_\gamma K_{n}(p,q) {Y}(q) \psi_j(q)  = \le\{\begin{array}{cc}
\psi_j(p) & j\leq n-1\\
0 & j\geq n
\end{array}\ri.
,\qquad
\int_\gamma K_{n}(p,q) {Y}(p) \psi_j^\vee(p)  = \le\{\begin{array}{cc}
\psi_j^\vee(q) & j\leq n-1\\
0 & j\geq n
\end{array}\ri.
\ee
 and the kernel is reproducing:
\be
\int_{q\in\gamma} K_n(p,q) {Y}(q) K_n(q,t) = K_n(p,t) \ ,\ \ \ \int K_{n}(p,p){Y}(p)=n.
\ee

In order to discuss a ``Christoffel--Darboux'' type of theorem one needs more structure; the nested sequence of divisors $\scr D_n$ will be taken with support within the pole divisor $\Xi$ of the projection map $Z$ of degree $r$. More precisely, 
\begin{assumption}
\label{assumeD}
 We consider the sequence of divisors $\scr D_n$ such that $\scr D_{kr+\ell} = k\Xi + \scr R_\ell$ with $\scr R_\ell$ positive divisor of degree $j\leq r-1$, $\scr  R_j \leq \Xi$.
\end{assumption}
\begin{example}
If $\Xi= \sum_{j=1}^r \infty^{(j)}$ is a simple divisor, then we may choose 
\be
\scr D_{kr+\ell} =  \sum_{j=1}^\ell (k+1) \infty^{(j)} + 
\sum_{j=\ell+1}^r k \infty^{(j)} , 
\ee
namely we increase  the multiplicities of the points at infinity sequentially and periodically.
Another example is $\Xi= (r-1)\infty^{(1)}+ \infty^{(2)}$; then we may choose 
\be
\scr D_{kr+\ell} = k\Xi + \ell \infty^{(1)} \ \ \ \  \ell \leq r-1.
\ee
A third example is simply $\Xi= r\infty$, namely, the map $Z$ has a single pole of order $r$; here the only choice we have is to take $\scr D_n = n\infty$.
\end{example}

Let the Assumption \ref{assumeD} prevail; then clearly $Z(p) \psi_n(p) \in \scr P_{n+r}$ and similarly $Z(p) \psi_n^\vee(p) \in \scr P^\vee_{n+r}$ just by inspection of the degree of the poles. 
Let us denote by $\vec \Psi$ the infinite column vector $(\psi_0,\psi_1,\dots)$ and $\vec \Psi^\vee$ the infinite row-vector $(\psi_0^\vee,\psi_1^\vee,\dots)$. Thanks to the fact that the multiplication by $Z$ is symmetric 
\be
\label{Zsymm}
\langle Z\varphi, \varphi^\vee\rangle = 
\langle \varphi,Z \varphi^\vee\rangle ,  \ \ \ \forall \varphi\in \scr P,  \ \ \forall \varphi^\vee \in \scr P^\vee,
\ee
it follows that the there is an $\mathbb N\times \mathbb N$ matrix $\mathbf Z$ such that 
\be
\label{multZ}
Z(p)\vec \Psi(p) = \mathbf Z \vec \Psi(p),\qquad Z(p)\vec \Psi^\vee(p) = \vec \Psi^\vee(p) \mathbf Z.
\ee
The  lemma hereafter follows  immediately from \eqref{Zsymm} and from the degree counting mentioned above.
\begin{lemma}
The $\mathbb N\times \mathbb N$ matrix $\mathbf Z$ is a finite-band matrix of size $2r+1$ with $\mathbf Z_{a,b}=0$ if $|a-b|> r$. 
\end{lemma}
It is convenient to view $\mathbf Z$ as a tridiagonal block matrix with blocks of size $r$ and similarly to segment the vectors $\vec \Psi, \vec \Psi^\vee$ into pieces of length $r$:
\be
\vec \Psi_k = [\psi_{kr}, \dots,\psi_{kr+r-1}]^t, \ \ \ \vec \Psi_k^\vee =  [\psi_{kr}^\vee, \dots,\psi_{kr+r-1}^\vee], \qquad k=0,1,\dots. 
\ee
Then \eqref{multZ} can be written 
\bea
Z(p)\vec \Psi_k(p) &= \AA_{k+1} \vec\Psi_{k+1}(p) + \BB_k \vec\Psi_k(p) + \CC_{k-1}\vec\Psi_{k-1}(p),
\\
Z(p)\vec \Psi_k^\vee(p)
& = \vec\Psi^\vee_{k-1}(p) \AA_{k} + \vec\Psi^\vee_k(p) \BB_k + \vec\Psi^\vee_{k+1}(p)\CC_{k},
\eea
where $\AA_k, \BB_k, \CC_k$ are $r\times r$ matrices.
Then we have the simple Corollary (whose proof is immediate)
\begin{corollary}[``CD identity"]
\label{CDI}
The CDS kernels for $n= \ell \,r$ can be written as follows
\be
K_{\ell r}(p,q) = \frac { \vec \Psi_{\ell}^\vee(q)\HH_\ell^{-1} \AA_{\ell+1} \vec \Psi_{\ell+1}(p) - \Psi_{\ell+1}^\vee(q) \CC_\ell \HH_\ell^{-1}\vec \Psi_{\ell}(p) }{Z(p)-Z(q)}
\ee
where $\HH_\ell = {\rm diag} (h_{\ell r}, \dots, h_{\ell r+r-1})$. 
\end{corollary}

\section{Matrix (bi)orthogonality}
\label{secMBOPs}
We consider a sequence of  divisors $\scr D_n$ chosen according to Assumption \ref{assumeD}.
The first important realization, contained in the next theorem, is that then the sections of $\scr P_n \setminus \scr P_{r-1}$ can be written as linear combinations of sections in $\scr P_{r-1}$ with polynomial coefficients in $Z$. 

We introduce the following notation; 
we choose a dissection, $\dot \C = \C \setminus \Sigma$, of the $z$--plane so that we can represent the Riemann surface $\mathcal C$ by $r$ copies of the $z$--plane appropriately identified along the various cuts of the dissection. 

For $z\in \dot \C$ we denote by $z^{(a)}, \ a=1,\dots, r$ the point in $Z^{-1} (z)\subset \mathcal C$ belonging to the corresponding sheet: two points $z^{(a)}, z^{(b)}$, $a\neq b$  coincide at  $z=c$ if and only if the value $c$ is a critical value of $Z$, namely $\d Z$ vanishes at some of the points in $Z^{-1}(c)$.
\bet
\label{Propideal}
Let $\{\varphi_j\}_{j=1}^r$ be any basis of holomorphic sections of $\mathcal X \otimes \sqrt{\mathcal K}(\Xi)$, with $\Xi$ the divisor of poles of the map $Z:\mathcal C\to \C\P^1$. Let
\be
\Phi(z) = \le[
\begin{array}{ccc}
\varphi_1(z^{(1)} )& \dots & \varphi_1(z^{(r)})
\\
\vdots &&\vdots\\
\varphi_r(z^{(1)}) & \dots & \varphi_r(z^{(r)})
\end{array}
\ri]
\label{Phi}
\ee
Then the determinant $\det \Phi(z)$ vanishes only at the critical points of minimal order in the following sense; for any collection $\psi_1,\dots, \psi_r$ of {\it meromorphic sections} of $\mathcal X\otimes \sqrt{\mathcal K}$ (and poles away from the ramification points of the map $Z$) the ratio
\be
F(z):= \frac {\det \bigg[\psi_j(z^{(a)})\bigg]_{j,a=1}^r}{\det \Phi(z)} 
\label{F}
\ee
is a rational function of $z$ without poles at the critical points (i.e. the determinant in the numerator is ``divisible'' by $\det \Phi$)
\eet
We could not find an elementary proof that does not use the theory of Theta functions and thus, inevitably, the proof we provide in Appendix \ref{proof} is somewhat technical. However, the reader may find the  genus $1$ example  in Section \ref{genus1example} sufficiently illuminating. 
\bp
\label{propZ}
Let $n\geq r$ and $\varphi\in \scr P_n\setminus \scr P_{r-1}$. Let  $\{\varphi_j\}_{1\leq j\leq r}$ be any  basis of $H^0(\mathcal X \otimes \sqrt{\mathcal K}(\Xi))$. Then the following two statements hold.\\
{\bf [1]} We have the expansion 
\be
\label{Zcomb}
\varphi(p) =  \sum_{j=1}^{r} f_j(z) \varphi_j(p),\ \ \ \ z= Z(p)
\ee
where $f_j$ are polynomials of degree $\leq \lfloor\frac n r \rfloor$ and at least one of them is of degree exactly $ \lfloor\frac n r \rfloor$.
\\
{\bf [2]} Let $\{\psi_\ell\}_{\ell \geq 0}$ be an  arbitrary sequence such that $\psi_\ell$ spans $\scr P_\ell \mod \scr P_{\ell -1}$ and arrange the  coefficients of the expansion \eqref{Zcomb} in matrix form as follows
\be
\le[\begin{array}{c}
\psi_{k r}(p)\\
\vdots\\
\psi_{kr+r-1}(p)
\end{array}
\ri] = F^{(k)}(z) \le[\begin{array}{c}
\psi_{0}(p)\\
\vdots\\
\psi_{r-1}(p)
\end{array}
\ri],\ \ \ z=Z(p).
\ee
Then $F^{(k)}(z)$ is a polynomial matrix of degree $k$ in the indeterminate $z$ and the leading coefficient matrix is invertible. 
\ep
{\bf Proof.} {\bf [1]} 
It is clear that any linear combination \eqref{Zcomb} belongs to  $\scr P_n$ so we need to show the converse.  We give here a rather pedestrian proof of this fact.
For $z \in \dot \C$ let $\vec \varphi(z):= \big( \varphi(z^{(1)}), \dots, \varphi(z^{(r)})\big)$ be the vector of its values at all the preimages of $z$. 
If $\sigma$ is a closed loop starting and ending at $z$ and avoiding the critical values of $Z$, then the analytic continuation of $\vec \varphi$ along $\sigma$ yields the same vector up to a quasi-permutation matrix $P_\sigma$, where the multipliers are determined by the  $\mathcal X$ multiplier system.

Consider then the matrix $\Phi$ in \eqref{Phi}: under analytic continuation along $\sigma$ the matrix $\Phi$ undergoes the same transformation:
\be
\Phi(z^\sigma) = \Phi(z) P_\sigma.
\ee
Moreover $\Phi$ is invertible on $\dot \C$ because of Thm. \ref{Propideal}.  Thus we can consider $\vec F(z):= \vec \varphi(z) \Phi^{-1}(z)$. It follows at once that this vector extends to a single valued function of $z$ on $\C$ minus the critical points of $Z$, where it can have at most poles so that $\vec F(z)$ is at most a rational function of $z$. We now show that it has no poles at the critical points,  and thus it is  a polynomial whose degree is easily determined to be as claimed. 

The $j$-th entry of $\vec F(z)$ can be written as a ratio of determinants:
\be
f_j(z) = \frac {\det \Phi_j}{\det \Phi}
\ee
where $\Phi_j$ is the matrix obtained by replacing the $j$-th row of $\Phi$ with the vector $\vec \varphi$ (Cramer's rule). Then the Theorem \ref{Propideal} states that the ratio is locally analytic also at the critical values. It follows from the Liouville theorem that $f_j$ are polynomials.  The degree follows from simple counting of the order of the poles. \\
{\bf [2]} As per Assumption \ref{assumeD} we have $\div(\psi_{kr+\ell-1}) \geq - k \Xi - \scr R_{\ell}$. To simplify the considerations we write  $\Xi = \infty^{(1)} + \dots +\infty^{(r)}$ (the points may be repeated), and $\scr R_\ell = \infty^{(1)} + \dots + \infty^{(\ell)}$, $\ell \leq r-1$ with $\scr R_0$ being the empty divisor.  
Then the sections $\psi_0,\dots, \psi_{r-1}$ form a basis of $H^0(\mathcal X\otimes \sqrt{\mathcal K}(\Xi))$ and we can apply the result at point {\bf [1]} of this proof. 
Let us denote by $F^{(k)}_k$ the coefficient matrix of $z^k$ of $F^{(k)}(z)$ of the polynomial matrix:
\be
F^{(k)}(z) = z^k F^{(k)}_k  + G^{(k)}(z).
\ee
where $G^k(z)$ is a polynomial matrix of degree $\leq k-1$.  
Now, the sections $\psi_{kr}, \dots, \psi_{kr+r-1}$ span $\scr P_{kr+r-1}\mod \scr P_{k(r-1)}$ and the $r$ sections given by the entries of  $G^{(k)}(z) [\psi_0,\dots, \psi_{r-1}]^t$  fall within $ \scr P_{k(r-1)}$ by simple degree counting. 
Thus the $r$ entries of 
\be
z^k F^{(k)}_k \le[
\begin{array}{c}
\psi_0\\
\vdots\\
\psi_{r-1}
\end{array}
\ri]
\ee
must span the $r$--dimensional space $\scr P_{kr+r-1}\mod \scr P_{k(r-1)}$ and hence $F^{(k)}_k$ must be invertible.
\QED
\begin{remark}
\label{remdegree}
If we consider an arbitrary increasing sequence of divisors $\scr D_n$ of degree $n$ supported at $\Xi$ but we waive the requirement that $\scr D_{kr+\ell} = k\Xi + \scr R_\ell$, then the Proposition \ref{propZ} still holds save for the statement about the degrees. In this more general case the sequence of degrees of the polynomials $f_j$ depends not just on the degree of $\scr D_n$ but on the choice of points. 
As a simple example of this phenomenon consider the map $Z:\C \P^1 \to \C \P^1$ given by $Z(t) = \frac 1 2(t - \frac 1 t)$. As basis of $\sqrt{\mathcal K}(\Xi)$ (there is no nontrivial flat bundle on $\C \P^1$) we can choose $\varphi_1 = \sqrt{d t}$ and $\varphi_2 = \frac {\sqrt{\d t}}{t}$.
If $\scr D_n= n\infty_t$, for example then the degrees of the polynomials in \eqref{Zcomb} are clearly $n$ and not $\lfloor n/2\rfloor$  because  the function $Z$ has a simple pole at infinity.
\end{remark}
Similar statements to those of Theorem \ref{Propideal} and Proposition \ref{propZ} hold for  sections of $\scr P_n^\vee$ with the obvious modifications.  
The following corollary is then immediate
\begin{corollary}
\label{corMBOPs}
Let $\psi_j,\psi_j^\vee$ be the biorthogonal sections for the pairing \eqref{pairing}. 
Denote by $\Psi_k(z), \Psi_k^\vee(z)$ the following matrices
\bea
\Psi_k(z) = \le[
\begin{array}{ccc}
\psi_{rk}(z^{(1)}) & \dots & \psi_{rk}(z^{(r)})\\
\vdots && \vdots\\
\psi_{rk+r-1}(z^{(1)}) & \dots & \psi_{rk+r-1}(z^{(r)})
\end{array}
\ri]\ ,\qquad 
\Psi_k^\vee (z) = \le[
\begin{array}{ccc}
\psi_{rk}^\vee(z^{(1)}) & \dots & \psi_{rk+r-1}^\vee(z^{(1)})\\
\vdots && \vdots\\
\psi_{rk}^\vee(z^{(r)}) & \dots & \psi_{rk+r-1}^\vee(z^{(r)})
\end{array}
\ri]
\eea
Then there exist matrix--valued polynomials $P_k(z), P_k^\vee(z)$ of degree $k$ such that 
\be
\label{PPvee}
\Psi_k(z)  = P_k(z) \Phi(z)\ ,\ \ \ \Psi_k^\vee (z) = \Phi^\vee(z) P_k^\vee(z),
\ee
where $\Phi,\Phi^\vee$ are the matrices constructed out of a pair of  arbitrary holomorphic bases $\mathfrak B = \{\varphi_\ell\}_{\ell=1}^r$, $\mathfrak B^\vee=\{\varphi_\ell^\vee\}_{\ell=1}^r$ of $\mathcal X\otimes \sqrt{\mathcal K}(\Xi), \mathcal X^\vee \otimes \sqrt{\mathcal K}(\Xi)$ (respectively) as follows
\be 
\Phi_{ab}(z) = \varphi_a(z^{(b)});\ \ \ \ \Phi^\vee_{ab} = \varphi_b^\vee (z^{(a)}).
\label{defPhigen}
\ee
\end{corollary}

The matrix polynomials $P_k, P_k^\vee$ will be our matrix-orthogonal polynomials under a few additional assumptions.

Let $\gamma_0\subset \C$ be an oriented contour in the plane and $ \gamma = Z^{-1}(\gamma_0)$: let ${Y}(z)$ as before, a smooth function on $ \gamma$ and define
\be
\Lambda(z) = {\rm diag} \bigg({Y}(z^{(1)}), \dots, {Y}(z^{(r)})\bigg).
\ee
Note that if $\gamma$ consists of several connected components, we are free to choose the function ${Y}$ on each of them, namely, we are not requiring that ${Y}$ extends outside of $ \gamma$ but simply that it is smooth on each of its connected components.
We also define the {\bf matrix weight} on $\gamma_0$ by the formula 
\be
\label{weight}
W(z)\d z:= \Phi(z) \Lambda(z) \Phi^\vee(z).
\ee
In general $\Lambda$ {\it is not} the diagonalization of the weight matrix because $\Phi^\vee$ is not the inverse of $\Phi$. This reflects the arbitrariness of choices of bases $\mathfrak B, \mathfrak B^\vee$ which affect $W$ by the left/right action of $\mathrm {GL}_n(\C)$.  
We then have the simple theorem 
\bet
\label{thmMBOPs}
Assume that $D_n\neq 0$ in \eqref{bimo} so that the biorthogonal sections $\psi_n, \psi_n^\vee$ span $\scr P_n\mod \scr P_{n-1}$ and $\scr P^\vee_n\mod \scr P^\vee_{n-1}$ respectively.
The matrix--valued polynomial sequences $\{P_k(z)\}_{k\geq 0}$, $\{P_k^\vee (z)\}_{k\geq 0}$  defined in \eqref{PPvee} are a biorthogonal sequence for the matrix weight $W(z)\d z$ \eqref{weight}:
\be
\int_{\gamma_0} P_{k}(z) W(z) P_\ell^\vee(z)\d z= \delta_{\ell k}\HH_\ell, 
\ee
where $\HH_\ell$ is the diagonal matrix 
\be
\HH_\ell = {\rm diag} \bigg(h_{\ell r
}, \dots, h_{\ell r+r-1}\bigg).
\ee
We can normalize the two sequences to be monic 
\be
Q_k(z) = (F_{k}^{(k)})^{-1}P_k(z), \qquad 
Q^\vee_k(z) = (F_{k}^{\vee (k)})^{-1}P^\vee_k(z), \qquad 
\ee
where the matrices $F_{k}^{(k)}$ and $F_{k}^{(k)}$ are the leading coefficients of $P_n, P_n^\vee$ (respectively), invertible by Prop. \ref{propZ}. In this case the ``norming'' matrices are:
\be
\int_\gamma Q_{k}(z) W(z) Q^\vee_\ell(z)\d z= \delta_{\ell k}\wt \HH_\ell.\ \ \ \wt \HH_\ell =(F_{k}^{(k)})^{-1}\HH_\ell (F_{k}^{\vee (k)})^{-1}.
\ee
\eet
\noindent {\bf Proof.}
By construction we have 
\be
P_{k}(z) W(z) P_\ell^\vee(z) = \Psi_k(z) \Lambda(z) \Psi_\ell^\vee(z). 
\ee
Consider the entry $(a,b)$ of the above matrix; 
\be
\le(\Psi_k(z) \Lambda(z) \Psi_\ell^\vee(z)\ri)_{a,b} = \sum_{j=1}^r 
\psi_{kr+a}(z^{(j)}) {Y}(z^{(j)}) \psi^\vee_{\ell r +b}(z^{(j)}).
\ee
Integrating this expression on $\gamma$ is the same as the integral over $\wh \gamma$ of $\psi_{kr+a}{Y}(p) \psi_{\ell r +b}^\vee(p)$, which gives  $\delta_{\ell k} \delta_{a,b} h_{\ell r+a}$.
\QED
\paragraph{$\Z_2$ characters and orthogonality: symmetric weight matrices.}
Let us consider again the case $\mathcal X:\pi_1\to \Z_2$. In this case we can choose the same bases $ \varphi_\ell= \varphi_\ell^\vee$ so that $\Phi^\vee = \Phi^t$. Then we have $\Psi_k^\vee = \Psi_k^t$ and $P_k (z) = P_k^t(z)$. The weight matrix $W$ is then a {\it symmetric} matrix. 
Under a few additional assumptions one can easily consider a setup where $W$ is actually positive definite. It suffices to consider curves $\mathcal C$ admitting an antiholomorphic involution that fixes $\wh \gamma$ and positive--valued weight function ${Y}$. We can only consider meromorphic functions $Z:\mathcal C\to \C \P^1$  that intertwine the anti-involution on $\mathcal C$ with some antiholomorphic involution on $\C\P^1$ (i.e. the conjugate  by an arbitrary M\"obius transformation of the ordinary complex conjugation).  

\subsection{Matrix Christoffel-Darboux from scalar Christoffel--Darboux--Szeg\"o}
As observed in the special case where $genus(\mathcal C)=0$ in \cite{Charlier} (see Thm. 1.16 ibidem) the matrix-Christoffel-Darboux (MCD) reproducing kernel can be expressed in terms of the scalar one. 
Specifically we recall the standard definition
\bd
The matrix Christoffel--Darboux kernel is 
\be
\mathbb K_\ell(z,w):= \sum_{j=0}^{\ell-1} P_j^\vee(z)\HH_j^{-1} P_j(w),  \ \ \ell =1,\dots. 
\ee
\ed
The kernel satisfies
\be
\label{MCD}
\int_\gamma \mathbb K_\ell(z,t)W(t) \mathbb K_\ell(t,w)\d t = \mathbb K_\ell(z,w),\ \ \ \int_{\gamma} \mathbb K_\ell (z,z)W(z)\d z = \ell\,\1_r.
\ee
It follows directly from the definitions of $P_j(z), P_j^\vee(z)$ that we have 
\begin{corollary}
The matrix Christoffel-Darboux kernel  $\mathbb K_\ell(z,w)$ in \eqref{MCD} and the (scalar) Christoffel--Darboux--Szeg\"o kernel $K_n$  in \eqref{CDS} are related by 
\be
\le[\Phi^\vee(z)\mathbb K_\ell(z,w) \Phi(w)\ri]_{ab} =K _{\ell r }(w^{(a)}, z^{(b)}).
\ee
\end{corollary}
\subsection{Change of multiplier system $\mathcal X$}
We briefly recall, following ideas explored in  \cite{Bertola:EffectiveIMRN}, of the effect of change of multiplier system $\mathcal X \to \wt{\mathcal X}$.

Under our assumptions that $h^0(\mathcal X \otimes \sqrt{\mathcal K})=0=h^0(\wt{\mathcal X} \otimes \sqrt{\mathcal K})$ (which, we remind, is equivalent to $\Theta(\mathcal X)\neq 0 \neq \Theta(\wt{\mathcal X}))$), there are sections, $\varphi_0$ and $\wt \varphi_0$, of the two bundles with a simple pole and $g$ zeros. We choose the pole to be at the same point $p_0$.  Then the ratio, $F_{\mathcal X, \wt {\mathcal X} } = \frac {\wt \varphi_0}{\varphi_0}$, is   a meromorphic section  of the flat bundle $\wt{\mathcal X} \otimes \mathcal X^\vee$  with $g$ zeros and $g$ poles. 
Namely a multivalued meromorphic function on $\mathcal C$ such that, under analytic continuation around a loop $\sigma\in \pi_1$, it satisfies
\be
F_{\mathcal X, \wt {\mathcal X} }(p^\sigma) =\frac{\wt {\mathcal X}(\sigma)} { {\mathcal X}(\sigma)} F_{\mathcal X, \wt {\mathcal X} }(p).
\ee
This function can be simply written in terms of the Riemann Theta function as 
\be
F_{\mathcal X, \wt {\mathcal X} }(p) ={\rm e}^{\sum_{j=1}^g (\wt \alpha_j-\alpha_j) \mathfrak A_j(p)} \frac {\Theta\le(p -p_0+
 \wt{\mathcal X}\ri)} {\Theta\le(p-p_0 + {\mathcal X}\ri)}
\ee
where we recall that we identify the multiplier system  $\mathcal X$ ($\wt {\mathcal X}$) with the point in the Jacobian $\mathbb J(\mathcal C)$ (of the same symbol) 
\be
\mathcal X = \bs \beta-\bs \tau \bs \a,\ \ , \ \ \mathcal X(a_j) = {\rm e}^{2i\pi \a_j}, \ \ \  \mathcal X(b_j) = {\rm e}^{2i\pi \b_j},
\ee
and similar formulas for $\wt {\mathcal X}$.
Let 
\be
\div \le(\frac {\wt \varphi_0}{\varphi_0}\ri) = \wt {\scr D} - \scr D,\ \ \ \deg \scr D = \deg{ \wt{\scr D}} = g.
\ee
Now, any holomorphic $\varphi$ section in $H^0(\mathcal X\otimes \sqrt{\mathcal K}(\Xi))$ can be written as $\varphi (p) = f(p) \varphi_0$ with $f(p)$ a meromorphic function with $\div f \geq -\Xi+ \scr D$.  Let us set
\be
D(z):= {\rm diag} \le( F_{\mathcal X, \wt{\mathcal X}}(z^{(1)}),\dots, F_{\mathcal X, \wt{\mathcal X}}(z^{(r)}) \ri) .
\ee
\bp
Given two flat bundles $\mathcal X, \wt {\mathcal X}$ and corresponding bases $\mathfrak B, \mathfrak B^\vee, \wt {\mathfrak B},\wt{ \mathfrak B}^\vee$  of the corresponding spaces, let $\Phi,\Phi^\vee$ be defined as in \eqref{defPhigen} 
and similar expressions defining the tilde counterparts.
For a given weight function ${Y}:\gamma\to \C$ the weight matrices $W(z)$, $\wt W(z)$ are related by 
\be
\label{WWt}
\wt W(z) = L(z) W(z) R(z),
\ee
with 
\be
L(z):= \wt \Phi(z) D^{-1}(z) \Phi^{-1}(z),\qquad R(z):= (\Phi^\vee(z))^{-1} D(z) \wt \Phi^\vee(z) .
\ee
Both $L(z), R(z)$ are rational functions with poles at $Z(\wt {\scr D})$ and $Z(\scr D)$, respectively.
\ep
{\bf Proof.}
Since both $D$ and $\Lambda$ are diagonal matrices we can write the weight matrix as 
\be
W(z) := \Phi(z) D(z) \Lambda(z) D^{-1}(z) \Phi^\vee(z).
\ee
 Then \eqref{WWt} immediately follows. The important information that the theorem is conveying is that both $L(z):=\wt \Phi(z) D^{-1}(z) \Phi^{-1}(z)$ and $R(z):=  (\Phi^\vee(z))^{-1} D(z) \wt \Phi^\vee(z) $ are rational (i.e. single--valued)  functions of $z$ with poles only at the $Z$--projection of the two divisors $\scr D, \wt {\scr D}$ of degree $g$.
We inspect $L(z)$ with the considerations extending to $R(z)$ by swapping the roles of $\mathcal X,\wt {\mathcal X}$. 
The matrix $\Phi(z) D(z)$ is of the form 
\be
[\Phi(z) D(z)]_{ab} =\varphi_a(z^{(b)}) \frac {\varphi_0^\vee(z^{(b)})} {\varphi_0(z^{(b)})},\ \ \ \ 
a,b=1,\dots, r.
\ee
For each $a=1,\dots, r$ the  $a$-th row is obtained by evaluating $\varpi_a(p):= \varphi_a(p) \frac {\wt \varphi_0(p)}{\varphi_0(p)}$ at the points $Z^{-1}(z)$. It is apparent that $\varpi_a$ is a {\it meromorphic} section of $\wt {\mathcal X}\otimes \sqrt{\mathcal K}(\Xi)$ with poles at $\scr D$ coming from the zeros of $\varphi_0$. Then $\Phi(z) D(z) \wt \Phi(z)^{-1}$ is single--valued by Theorem \ref{Propideal} (applied to the basis $\wt {\mathfrak B}$). 
\QED

\subsection{Abelianization: from matrix to scalar orthogonality}
\label{abelianization}
The above construction starts with a Riemann surface and a projection $Z$. It is interesting to give a general framework where we start from the weight matrix $W(z)$ itself. 

This can be done at least  under the  assumption that $W(z)$ is a rational matrix; in this case there exists a well established general spectral theory which we briefly summarize below. In the context of MOPs, this precise approach has been investigated recently in \cite{Charlier} and hence we will not delve into all details.

The reader that would like to have a geometrical approach to the so--called ``(inverse) spectral problem" (such is the issue at hand) should consult for example Chapter 5 of \cite{BabelonBook} for a recent review, or the papers \cite{Adams1, Adams2, Bertola:EffectiveIMRN} 

We preface this construction by warning that it is not quite the inverse construction to the previous one; the reason is that we will use as weight on the Riemann surface the eigenvalues of $W$ while in the previous section the matrix $\Lambda$ is not the matrix of eigenvalues of $W$ in \eqref{weight} in general.  \\[1pt]

The spectral curve $\mathcal C$ is simply the locus of the characteristic polynomial of $W$:
\be
\mathcal C:=\Big  \{(y,z)\in \C^2: \ \ \  \det (y\1_r - W(z)) =0\Big\}.
\ee
The function ${Y}$ used in the previous section is now the eigenvalues of $W$ on the spectral curve, and $Z$ is the natural projection to the $z$--plane.  
The curve $\mathcal C$ can be compactified under general assumptions and we denote by $g$ the genus of this compactification. 

If $(y,z)\in \mathcal C$ then the matrix $M(y,z):= y\1_r -W(z)$ is singular and generically its {\it adjugate} matrix $\wt M(y,z)$  (the transposed matrix of cofactors) has rank one. Any nonzero column of $\wt M$ is then a right eigenvector and any nonzero row is a left eigenvector. Note that we are not assuming that $W$ be a symmetric matrix. 

Choose one such row $\vec R$ and  column $\vec C$ : choose one of the entries of each (say, the first) and  then it can be seen that they have $g+r-1$ zeros. The divisors of these zeros determine  line bundles of the same degree $\scr L, \scr L^\vee$ dual to each other.  Consider the case of $\scr L$; the divisor is linearly equivalent to a divisor of the form $\scr D-\infty^{(1)} + \Xi$, where $\scr D$ is positive and of degree $g$ and $\infty^{(1)}$ is one of the preimages of $z=\infty$.  
Let us fix a normalized basis of holomorphic differentials on $\mathcal C$ and the corresponding Abel map $\mathfrak A$ with basepoint $\infty^{(1)}$. We denote by ${\bs \alpha, \bs \beta}$ the (real) characteristics of $\scr D$ as follows, 
\be
\mathfrak A(\scr D) = {\bs \beta} - \bs \tau {\bs \alpha},\ \ \ {\bs \tau}_{ij} = \oint_{b_i} \omega_j  = \oint_{b_j} \omega_i, \ \ \ {\bs \beta,\bs \alpha} \in \R^g.
\ee
Let us denote by $\mathcal X = \mathfrak A(\scr D) = {\bs \beta} - \bs \tau {\bs \alpha}$; this point defines also a multiplier system which we denote by the same symbol $\mathcal X:\pi_1(\mathcal C)\to U(1)$ with slight abuse of notation according to the formula:
\be
\mathcal X(a_\ell) = {\rm e}^{2i\pi \alpha_\ell},\ \ \ \mathcal X(b_\ell)= {\rm e}^{2i\pi \beta_\ell},  \ \ \ell=1,\dots, g.
\label{318}
\ee
For simplicity let us assume at first  that $\Xi = Z^{-1}(\infty) = \sum_{j=1}^r \infty^{(j)}$ is a simple divisor; if $S(p,q):= S_{\bs \alpha,\bs \beta}(p,q)$ denotes the Szeg\"o\ kernel with the indicated characteristics (see Appendix \ref{theta}) then we can follow \cite{Korotkin} and define the matrix 
\be
\Phi_{j,a} (z) = \frac {S( z^{(a)}, \infty^{(j)})}{\sqrt{\d \zeta_j}} ,\ \ \ \Phi^\vee_{b,k}:=\frac{ S( \infty^{(k)},z^{(b)}) }{\sqrt{\d \zeta_k}}
\ee
where $\d \zeta_j$ denotes any local coordinate near $\infty^{(j)}$ in which to trivialize the semicanonical bundle $\sqrt{\mathcal K}$. If we choose $\zeta = \frac  1 z $ (for simple divisor $\Xi$) then they satisfy 
\be
\Phi(z) \Phi^\vee(z) = \1_r \d z.
\label{phiinv}
\ee
To see this (see \cite{Korotkin}) one computes directly the $j,k$ entry:
\be
(\Phi(z) \Phi^\vee(z))_{jk} = \sum_{a=1}^r  
\frac {S( \infty^{(k)}, z^{(a)} )}{\sqrt{\d \zeta_k}}\frac{ S(z^{(a)}, \infty^{(j)})}{\sqrt{\d \zeta_j}} 
\ee
Since the expression is clearly invariant under permutation of the sheets, and single valued,  it must be a differential on the target curve of the $z$--plane; as such it can only have poles at $\infty$ due to the simple pole along the diagonal of the Szeg\"o kernel; however if $j\neq k$ the pole is  of order at most  $1$ and  since there are no nontrivial differentials with only one simple pole on the Riemann sphere, we conclude that  it vanishes.  For $j=k$ the expression has a double pole at $z=\infty$ due to the fact that the Szeg\"o\ kernel has a simple pole on the diagonal and hence it is proportional to $\d z$. To show that the proportionality is $1$ we recall that if $p,q$ are in the same coordinate patch $\zeta$ then 
\eqref{Szegonorm} applies.  Then, near $\infty_j$ the coordinate is $\zeta_j = \frac 1 z$ and hence one finds that the coefficient is  $1$ on all sheets. 

\paragraph{More general divisor $\Xi$.}
If $\Xi$ is not simple, then (see also App. \ref{theta}) we proceed as follows to construct sections of $\mathcal X\otimes \sqrt{\mathcal K}(\Xi)$. Suppose that $\infty^{(j)}$ has multiplicity $\nu_j$: then take local coordinate  $\zeta_j = z^{-\nu_j}$ and consider the expressions 
\be
\varphi(p,\zeta_j):= \frac{S(p,\zeta_j^{-1}(\zeta_j))}{\sqrt{\d \zeta_j}};\ \ \ \varphi^\vee(q,\zeta_j):= \frac{S(\zeta_j^{-1}(\zeta_j),q)}{\sqrt{\d \zeta_j}}.
\ee
Namely, we evaluate the Szeg\"o\ kernel in the coordinate $\zeta_j$ with respect the second variable (first, respectively). Consider the expression $\varphi(p,w)$; it is a section of $\mathcal X\otimes \sqrt{\mathcal K}$ in the variable $p$ with a simple pole at $\zeta_j^{-1}(w)$. Thus, it  depends parametrically on the complex variable $w$ in a neighbourhood of $w=\infty$. Similar considerations apply to  $\varphi(q,w)$ as a section of $\mathcal X^\vee\otimes \sqrt{\mathcal K}$.

 The dependence on $z$ is that of a Taylor series in $\zeta_j:= w^{-\frac 1 {\nu_j}}$ (i.e. a Puiseux Laurent series).
Then one defines 
\be
\label{baseses}
\varphi_{\ell;j}(p):=\frac 1{(\ell-1)!} \le(\frac {\d}{\d \zeta_j}\ri)^{\ell-1}\!\!\!\!\!\! \varphi(p,\zeta_j) \bigg|_{\zeta_j=0},\ \ 
\varphi_{\ell;j}(p)^\vee :=\frac {(-1)^{\ell} }{(\ell-1)!} \le(\frac {\d}{\d \zeta_j}\ri)^{\ell-1}\!\!\!\!\!\! \varphi^\vee(p,\zeta_j) \bigg|_{\zeta_j=0}
\ \ell= 1,\dots, \nu_j.
\ee
It follows from the properties of the Szeg\"o\ kernel that these are holomorphic sections of $\mathcal X\otimes \sqrt{\mathcal K}(\ell\infty_j)$  (and hence also holomorphic sections of $\mathcal X\otimes \sqrt{\mathcal K}(\Xi)$) with a single  pole of order $\ell$ at $\infty_j$ and normalized as 
\be
\varphi_{\ell;j}(p) = \zeta_j(p)^{-\ell} \le(1+ \mathcal O(\zeta_j(p))\ri) \sqrt{\d \zeta_j(p)}.
\ee
In particular they are clearly linearly independent. Similar considerations apply to $\varphi^\vee$.  Thus, for a general divisor $\Xi= \sum_{j=1}^k \nu_j \infty^{(j)}$  a pair of bases of holomorphic sections of  $\mathcal X\otimes \sqrt{\mathcal K}(\Xi)$  and $\mathcal X^\vee\otimes \sqrt{\mathcal K}(\Xi)$ are
\bea
&\mathfrak B:= \bigg\{\varphi_{\ell;j}(p), \ j=1,\dots, k;\ \ \ell = 1,\dots, \nu_j \bigg\}
\ 
&\mathfrak B^\vee := \bigg\{\varphi^\vee_{\ell;j}(p), \ j=1,\dots, k;\ \ \ell = 1,\dots, \nu_j \bigg\}.
\eea

We arrange the bases in the natural order implied by the indexing and construct the matrices
\be
\Phi(z):= \le[
\begin{array}{ccccc}
\varphi_{1;1}(z^{(1)}) & \dots & \varphi_{1;1}(z^{(r)}) \\
\vdots && \vdots\\
\varphi_{1;\nu_1}(z^{(1)}) & \dots & \varphi_{1;\nu_1}(z^{(r)}) \\
\hline\\
\vdots\\
\hline
\varphi_{k;1}(z^{(1)}) & \dots & \varphi_{k;1}(z^{(r)}) \\
\vdots && \vdots\\
\varphi_{k;\nu_k}(z^{(1)}) & \dots & \varphi_{k;\nu_k}(z^{(r)}) \\
\end{array}
\ri]
\\
\Phi^\vee(z):= \le[
\begin{array}{c|c|c}
\varphi^\vee_{1;1}(z^{(1)})  \dots  \varphi^\vee_{1;\nu_1}(z^{(1)}) &\dots &\varphi^\vee_{k;1}(z^{(1)})  \dots  \varphi^\vee_{k;\nu_k}(z^{(1)}) \\
\vdots && \vdots\\
\varphi^\vee_{1;1}(z^{(r)})  \dots  \varphi^\vee_{1;\nu_1}(z^{(r)}) &\dots &\varphi^\vee_{k;1}(z^{(r)})  \dots  \varphi^\vee_{k;\nu_k}(z^{(r)}) 
\end{array}
\ri]
\label{defPhigen2}
\ee

Then we have 
\bp
\label{propPhiPhivee}
The product $\Phi(z)\Phi^\vee(z)$ is a constant matrix times $\d z$ of the form:
\be
\Phi(z)\Phi^\vee(z) = {\rm diag} \bigg( J_{\nu_1},\dots,J_{\nu_k}\bigg)\d z = J \d z
\ee
where each block $J_\ell$ is  a constant (in $z$) antitriangular matrix with unit on the antidiagonal entries: $(J_{\ell})_{a, r+1-b}=-1$ and nonzero entries only for $a+b\geq r+1$.
\ep
{\bf Proof.} Without great loss of generality and to avoid a plethora of indices, we consider the case $\Xi= r\infty$ (i.e. with $k=1$ and $\nu_1=r$). The basis $\varphi_j$ defined in \eqref{baseses} behaves as (we write it directly in the coordinate form)
\be
\varphi_j(z) =\zeta^{-j}\big(1 + \mathcal O(\zeta)\big) {\sqrt{\d\zeta}},\ \ \  j=1,\dots r
\ee
and $\varphi^\vee_j$  has the same behaviour. Here $\zeta = z^{-\frac 1 r}$ is the local coordinate near $\infty\in \mathcal C$.  On the different sheets the different behaviours are obtained by multiplying $\zeta$ by the $r$ roots of unity. Thus 
\be
\sum_{a=1}^r \varphi_j(z^{(a)}) \varphi_k(z^{(a)}) = \zeta^{-j-k} \le(\sum_{a=1}^{r} 
{\rm e}^{\frac {2i\pi a} r (j+k-1)}   + c_1   {\rm e}^{\frac {2i\pi a} r (j+k-2)}\zeta +\mathcal O(\zeta^2)\ri)\d \zeta. 
\ee
We know, a priori, that the result is a series in integer powers of $z$ (this can also be seen by the fact that the sum of the roots of unity is zero). Moreover we know that the result must be a polynomial times $\d z$ because there are no singularities in the finite $z$-plane. Thus we conclude that we can have a nonzero result only for $2r-1\geq j+k-1 \geq r$. For $j+k-1=r$ we get $r \frac {\d \zeta}{\zeta^{r+1}} = -\d z$.
\QED

\begin{remark}
\label{remphinv}
It is clear that by a triangular change of basis we can arrange for basis sections 
in such a way that $\Phi\Phi^\vee = \1 \d z$.
\end{remark}
We now assume that the bases have been chosen so that the Remark \ref{remphinv} applies and hence $\Phi\Phi^\vee = \1 \d z$, $\Phi^{-1} = \frac 1 {\d z} \Phi^\vee$. 
In this case the matrix $\Lambda$ is indeed the diagonal matrix of eigenvalues of $W =  \Phi(z) \Lambda(z) \Phi^\vee(z) $ (while in the previous construction it needed not be). 
 
The construction of the matrix biorthogonal polynomials then proceeds by construction  the (scalar)  orthogonal sections as in Section \ref{BOPs} and then the resulting matrix biorthogonal polynomials as in  Section \ref{secMBOPs}.

\section{Examples: genus zero and one}
\label{secExamples}
The reader without familiarity with the above algebro-geometric language  may find respite in contemplating the following two classes of examples in genus zero and one.
\subsection {Genus zero cover}
We assume that $\mathcal C= \C \P^1$, with affine coordinate $t$ and $Z(t)$ is simply a rational function of degree $r$.

We choose a path $\gamma_0$ in the $z$-plane and its pullback $ \gamma = Z^{-1}(\gamma_0)$. 
Given any smooth function ${Y}(t)$ defined on $\gamma$ one can orthogonalize the spaces $\scr P_n, \scr P_n^\vee$ and produce instances of matrix orthogonal polynomials. 

To keep the example concrete and simple, we take $Z(t)$ to be a polynomial of degree $r$: 
\be
Z(t) = t^r  + a_1 t^{r-1}+ \dots + a_r.
\ee
Then $\Xi= r\infty_t$ (we use the subscript to indicate that this is the point at infinity in the $t$--plane).
Since the map $Z$ has a single pole, we necessarily must set $\scr D_n = n\infty_t$ so that $\scr P_n$ is identified with polynomials of degree $n-1$ in the $t$--variable.

Then the matrix weight on the $z$ plane is given by 
\be
W(z) \d z = \Phi(z) \Lambda(z) \Phi^\vee(z), \ \ \ \ z\in \gamma_0,
\label{weightt}
\ee
where $\Phi, \Phi^\vee$ are constructed as in \eqref{defPhigen2} and   
\be
\Lambda(z) := {\rm diag} \big({Y}(z^{(1)}), \dots, {Y}(z^{(r)})\big).
\ee
Note that in general  $\Lambda$ is not the matrix of eigenvalues of $W$. 

Now the Proposition \ref{propZ} translates to a simple construction in commutative algebra; a polynomial $P(t)$ of degree $n\geq r$ can be written as 
\be
P(t) = Q_0(t)Z(t) + R_0(t),\ \ \ \deg R_0(t)\leq r-1,  \ \ \deg Q_0= n-r.
\ee
by long division. If $n-r\geq r$ then we can perform the long division on $Q(t)$ and obtain 
\be
P(t) = Q_1 (t)Z^2(t) + R_1 (t)Z(t) + R_0(t),\ \ \ \deg Q_1 = n-2r, \ \ \deg R_1 \leq r-1.
\ee 
We clearly can repeat this procedure until we reach a step $k= \lfloor \frac n r\rfloor$ where $\deg Q_k\leq r-1$; thus we can uniquely write
\be
P(t) = \sum_{j=0}^k R_j(t) Z^j(t),\ \ \ \deg R_j\leq r-1.
\ee
We can now express each of the polynomials $R_j(t)$ as a linear combination with {\it constant} coefficients of the polynomials $p_0,\dots, p_{r-1}$ and rearrange the sum so as to get an expression of the form 
\be
P(t)\sqrt{d t} = \sum_{\ell=1}^{r} f_\ell(Z) \psi_{\ell-1}(t),\ \ \ \psi_{\ell}(t) = p_{\ell}(t) \sqrt{\d t}.
\ee
This is the elementary version of the Proposition \ref{propZ}. Thus we are implicitly choosing the bases $\varphi_\ell = \psi_{\ell-1}$, $\ell =1,\dots, r$. 

Now consider the orthogonal polynomial $p_{nr+a}$ with $0\leq a \leq r-1$; by the above considerations we can write it as 
\be
\psi_{nr+a-1}(t)=
p_{nr+a-1}(t)\sqrt {\d t} = \sum_{\ell=1}^r P_{n;a,\ell} (z) \varphi_{\ell}(t),\ \ \ z=Z(t),
\ee
and thus 
\be
\Psi_n(z) = \Big[\psi_{nr+a-1}(z^{(\ell)}) \Big]_{a,\ell=1}^r= P_n(z) \Psi_0(z),  \ \ \ \ \ 
\Psi_n^\vee(z) = \Psi_n(z)^t.
\ee
We thus have 
\bp
The matrix polynomials $P_n(z)$ are orthogonal with respect to the matrix weight $W(z)\d z$ given in \eqref{weightt}:
\be
\int_{\gamma_0} P_n(z) W(z) P^t_{m}(z) \d z = \delta_{nm} H_n\ ,\ \ \ \ 
H_n= {\rm diag} \big (h_{nr}, \dots, h_{nr+r-1}\big).
\ee
\ep

\subsubsection{Laguerre matrix polynomials}
\label{secLaguerre}
We give a simple example of the above scheme, but clearly one can devise infinitely many such simple examples.
Consider the map $Z(t)= (t-c)^2$, $c\leq 0$. We can invert the map $Z$ as 
\be
t_{1,2}(z) =  c\pm \sqrt{z} ,\ \ \ \d t =  \frac {\pm \d z}{ {2} \sqrt{z}}
\ee
We consider the Laguerre polynomials  $\mathrm L_j^\alpha$  on $R_+\subset \C \P^1_t$:
\bea
\mathrm L_n^\alpha(t) &=\frac{t^{-\alpha} {\rm e}^t}{n!}  \frac {\d^n}{\d t^n} {\rm e}^{-t} t^{n+\alpha}\\
&\int_{\R_+} \mathrm L_n^\alpha(t) \mathrm L_m^\alpha(t) {\rm e}^{-t}t^\alpha \d t = \delta_{nm} \frac {\Gamma (n+\alpha+1)}{n!}.
\eea
We choose as basis for $\sqrt{\mathcal K}(\Xi)$ the sections $\psi_0(t) = \sqrt{\d t}$ and $\psi_1(t) =  \mathrm L_1^\a(t) \sqrt{\d t} =(1+\a-t)\sqrt{\d t}$, which we need to evaluate at $t_{1,2}$.
Then  we  find 
\bea
\Phi(z) =  \le[
\begin{array}{cc}
1 & i \\
\a+1-c-\sqrt{z}  & i(\a+1-c+\sqrt{z})
\end{array}
\ri]\frac {\sqrt{\d z}} {\sqrt{2} z^\frac 1 4},\ \ \ \Phi^\vee(z) = \Phi^t(z).
\eea
We note that 
\be
\Phi(z) \Phi^\vee(z) = \le[\begin{array}{cc}
0 & -1\\ -1 & 2(c-1-\a)
\end{array}\ri]\d z
\ee
The half-line is mapped to $[c^2,\infty)\subset \C\P^1_z$.  We then set  
\be
\Lambda(z):= 
\le[
\begin{array}{cc}
t_1(z)^{\alpha} {\rm e} ^{-t_1(z) } & 0\\ 0& 0
\end{array}
\ri]
\label{Lag}
\ee
Note that since $c<0$ we have $t_1(c^2)= c+ \sqrt{c^2} = 0$. 
A short computation yields the weight matrix: 
\be
W_L(z) \d z= \Phi (z) \Lambda(z) \Phi^\vee(z)=
\le[\begin{array}{cc}
1 &\a+1-c-\sqrt{z}\\
\a+1-c-\sqrt{z}  &(\a+1-c-\sqrt{z})^2 
\end{array}\ri]
(c+\sqrt{z})^\alpha \frac{{\rm e}^{-c-\sqrt z}}{2\sqrt z} \d z.
\ee
Finally we express an arbitrary polynomial $p(t)$ in terms of a linear combination of $\mathrm L_0^\a(t)=1$ and $L_1^\a(t)=\a+1-t$ with coefficients being polynomials in $z=Z(t)$.
A direct verification shows 
\be
p(t) = \res{s=\infty}\frac{ p(s) (s+1+\a-2c)\d s}{z-Z(s)} +(1+\a-t) \res{s=\infty} \frac {p(s)\d s}{(Z(s)-z)} 
\ee
which allows to explicitly compute the MOPs:
\be
P_j(z) = \le[
\begin{array}{cc}
\res{s=\infty}\frac{ \mathrm L^\alpha_{2j}(s) (s+1+\a-2c)\d s}{z-Z(s)} & \res{s=\infty} \frac {\mathrm L^\alpha_{2j}(s)\d s}{Z(s)-z}\\
\res{s=\infty}\frac{ \mathrm L^\alpha_{2j+1}(s) (s+1+\a-2c)\d s}{z-Z(s)} & \res{s=\infty} \frac {\mathrm L^\alpha_{2j+1}(s)\d s}{Z(s)-z}\\
\end{array}.
\ri]
\ee
By construction these polynomials satisfy 
\be
\int_{c^2}^\infty P_j(z) W_L(z) P_k^t(z) \d z = \delta_{jk} \le[
\begin{array}{cc}
\frac {\Gamma(2j+\alpha+1)}{(2j)!} & 0\\
0& 
\frac {\Gamma(2j+\alpha+2)}{(2j+1)!}
\end{array}\ri].
\ee

Note that $\Lambda$ in \eqref{Lag} is {\it not} the matrix of eigenvalues for $W_L$.
This weight matrix (as well as the one in the next example) does not fall within the purview of \cite{Charlier} because it is not a rational matrix.

\subsubsection{Hermite matrix polynomials}
\label{secHermite}
To give another  completely explicit example let us consider the Hermite polynomials in the $t$--plane and the same degree $2$ map $Z(t) = (t-c)^2$, $c\in \R$:
\be
{\mathrm h}_n(t) = (-1)^n e^{t^2} \frac {\d^n}{\d t^n} {\rm e}^{-t^2},\ \ {\mathrm h}_0(t)=1,  \ \ {\mathrm h}_1(t)=2t,\dots
\ee
We choose as basis for $\sqrt{\mathcal K}(\Xi)$ the sections $\psi_0(t) = \sqrt{\d t}$ and $\psi_1(t) = \mathrm h_1(t) \sqrt{\d t} =2t\sqrt{\d t}$, which we need to evaluate at $t_{1,2}$.
 We then find
\be
\label{psi_0}
\Phi(z) = \le[
\begin{array}{cc}
1 & i \\
2c+2\sqrt{z}  & 2i(c-\sqrt{z})
\end{array}
\ri]\frac {\sqrt{\d z}} {\sqrt{2} z^\frac 1 4}
\ee
We note that 
\be
\Phi(z) \Phi(z)^t = \le[\begin{array}{cc}
0 & 2\\ 2 & 8c
\end{array}\ri]\d z
\ee
Note that the matrix in the right side can be factored as $M M^t$ with 
\be
M = \le[
\begin{array}{cc}
\frac{-i}{\sqrt{2 c}} & \frac 1{\sqrt{2c}}\\ 0 & \sqrt{8c}
\end{array}
\ri] ,\ \ \  M M^t = \le[\begin{array}{cc}
0 & 2\\ 2 & 8c
\end{array}\ri]
\ee
This would allow us to re-define $\Phi \mapsto = M^{-1} \Phi$ and $\Phi^\vee \mapsto  \Phi^\vee M^{-T}$ in such a way that actually $\Phi \Phi^\vee = \1 \d z$. However the final weight matrix would not be real--valued. 

The image of the oriented path $ \gamma =\R\subset \C \P^1_t$ under $Z$ is the half-line $\gamma_0 = [c^2,\infty)\subset \C \P^1_z$ covered twice in opposite directions. 
We then set  
\be
\Lambda(z):= 
\le[
\begin{array}{cc}
{\rm e} ^{-(c + \sqrt{z})^2 } & 0\\ 0& -{\rm e}^{-(c-\sqrt{z})^2}
\end{array}
\ri]
=
\le[
\begin{array}{cc}
{\rm e} ^{-(c^2+z  +2c \sqrt{z}) } & 0\\ 0&- {\rm e}^{-(c^2+z-2c\sqrt{z})}
\end{array}
\ri]
\ee
The minus sign is due to the fact that the orientation is flipped on the second sheet.
A short computation using $\Phi$ as in \eqref{psi_0} yields the weight matrix: 
\be
W(z) \d z= \Phi (z) \Lambda(z) \Phi^t(z)=
\le( \frac{\cosh(2 c \sqrt{z})} {\sqrt{z}}
\le[\begin{array}{cc}
1 & 2c\\ 2c & 4c^2+4z
\end{array}\ri]
- 2{\sinh(2c\sqrt{z})}
\le[\begin{array}{cc}
0& 1\\ 1 & 4c
\end{array}\ri]
\ri) {\rm e}^{-z-c^2} \d z.
\ee
Note that for $c=0$ the matrix $W$ is diagonal and also the matrix orthogonal polynomials become diagonal. This is simply a manifestation that in this case $z=t^2$ and the Hermite polynomials respect the parity. 

Finally we express an arbitrary polynomial $p(t)$ in terms of a linear combination of $\mathrm h_0(t)=1$ and $\mathrm h_1(t)=2t$ with coefficients being polynomials in $z=Z(t)$.
A direct verification shows 
\be
p(t) = \mathrm h_0(t)\res{s=\infty}\frac{ p(s) (s-2c)\d s}{z-Z(s)} + \mathrm h_1(t) \res{s=\infty} \frac {p(s)\d s}{2(z-Z(s))}
\ee
which allows to explicitly compute the MOPs;
\be
P_j(z) = \le[
\begin{array}{cc}
\res{s=\infty}\frac{ \mathrm h_{2j}(s) (s-2c)\d s}{z-Z(s)} & \res{s=\infty} \frac {\mathrm \mathrm h_{2j}(s)\d s}{2(z-Z(s))}\\
\res{s=\infty}\frac{ \mathrm h_{2j+1}(s) (s-2c)\d s}{z-Z(s)} & \res{s=\infty} \frac {\mathrm h_{2j+1}(s)\d s}{2(z-Z(s))}\\
\end{array}
\ri].
\ee
The first few are $P_0=\1$ and 
 \bea
P_1&= \left[ \begin {array}{cc} -4\,{c}^{2}+4\,z-2&4\,c
\\ \noalign{\medskip}-16\,{c}^{3}+16\,cz&12\,{c}^{2}+4\,z-6
\end {array} \right]
\\
P_2&=  \left[ \begin {array}{cc} 16\,{z}^{2}+ \left( 32\,{c}^{2}-48 \right) 
z-48\,{c}^{4}+48\,{c}^{2}+12&32\,{c}^{3}+32\,cz-48\,c
\\ \noalign{\medskip}-128\,{c}^{5}+320\,{c}^{3}+128\,c{z}^{2}-320\,cz&
16\,{z}^{2}+ \left( 160\,{c}^{2}-80 \right) z+80\,{c}^{4}-240\,{c}^{2}
+60\end {array} \right] 
\eea
{\tiny
\bea
P_3&= \left[ \begin {array}{cc} 64  {z}^{3} \!+\!  \left( 576  {c}^{2} \!-\! 480
 \right) {z}^{2} \!+\!  \left(  \!-\! 320  {c}^{4} \!-\! 960  {c}^{2} \!+\! 720 \right) z \!-\! 320
  {c}^{6} \!+\! 1440  {c}^{4} \!-\! 720  {c}^{2} \!-\! 120&192  c{z}^{2} \!+\!  \left( 640  {c
}^{3} \!-\! 960  c \right) z \!+\! 192  {c}^{5} \!-\! 960  {c}^{3} \!+\! 720  c
\\ \noalign{\medskip}768  c{z}^{3} \!+\!  \left( 1792  {c}^{3} \!-\! 5376  c
 \right) {z}^{2} \!+\!  \left(  \!-\! 1792  {c}^{5} \!+\! 6720  c \right) z \!-\! 768  {c}^{7}
 \!+\! 5376  {c}^{5} \!-\! 6720  {c}^{3}&64  {z}^{3} \!+\!  \left( 1344  {c}^{2} \!-\! 672
 \right) {z}^{2} \!+\!  \left( 2240  {c}^{4} \!-\! 6720  {c}^{2} \!+\! 1680 \right) z \!+\! 
448  {c}^{6} \!-\! 3360  {c}^{4} \!+\! 5040  {c}^{2} \!-\! 840\end {array} \right].
\eea
}
By construction these polynomials satisfy 
\be
\int_{c^2}^\infty P_j(z) W(z) P_k^t (z) \d z =\sqrt{\pi} \delta_{jk} \le[
\begin{array}{cc}
2^{2j} (2j)! & 0\\
0& 2^{2j+1}(2j+1)!
\end{array}\ri].
\ee

\subsection{Elliptic covers}
\label{genus1example}
\subsubsection{Genus 1}
%
%

%

%%%%%%%%%%%%%%%%%%%%

Let $\mathcal C$ be an elliptic curve, which we represent as the quotient $v\in \C / \Z + \tau \Z$ with $\tau = \frac {\omega_2}{\omega_1}\in \mathbb H_+$.
The curve can also be represented in Weierstra\ss\ form:
\be
y^2 = 4 z^3 - g_1 z-g_3 = 4(z-e_1)(z-e_2)(z-e_3).\label{ellyz}
\ee
The uniformization is given by the Weierstra\ss's $\wp$ function as $y =(2\omega_1)^{-3} \wp'(v), \ z=(2\omega_1)^{-2}\wp(v)$ and 
\bea
\omega_1 &= \int_{\infty}^{e_1} \frac {\d z}{y}, \ \ 
\omega_2 = \int_{e_1}^{e_2} \frac {\d z}{y},\ \ \ \tau:= \frac{\omega_2}{\omega_1}\\
\wp(v)&:= \frac 1 {v^2} + \sum_{m^2+n^2\neq 0\atop m,n\in \Z} \le(
\frac 1{(v+m+n\tau)^2} - \frac 1{(m+n\tau)^2}\ri)
\eea
and then the curve \eqref{ellyz} is identified with the $2$--torus (i.e. the {\it Jacobian} $\mathbb J$)  $v\in \C \mod \Z + \tau \Z$, with the point $z=\infty$ corresponding to $v=0$ and the three points $e_1,e_2,e_3$ corresponding with the half-periods $v=\frac 1 2, \frac \tau 2, \frac {1+\tau}2$, respectively.
The point at infinity in the $z$--plane is a double point corresponding to $v=0$.  

Since the function $z = Z(p) = \wp(v)$ has a single (double) pole at $v=0$, the only choice of sequence of divisors is $\scr D_n = n\infty$. 
Let $\theta_1$ be the Jacobi (odd) theta function 
\bea
\label{Jacobitheta}
\theta_1(v;\tau) ={\rm e}^{\frac{i\pi}4 \tau + i \pi (v+\frac 1 2)} \sum_{n\in \Z} {\rm e}^{i\pi n^2\tau + 2i\pi n (v+\frac 1 2 + \frac \tau 2)}.
\eea
The function $\theta_1$ is odd, vanishes at $v=n+m\tau$ of simple order and satisfies
\be
\theta_1(v+1)= -\theta_1(v), \ \ \theta_1(v+\tau) = -{\rm e}^{-2i\pi v-i\pi \tau } \theta_1(v). 
\ee

\paragraph{Szeg\"o kernel.}  

The Szeg\"o\ kernel  for a multiplier system $\mathcal X$ is expressed using the Jacobi theta function \eqref{Jacobitheta} as follows:
\be
S(v,w) = {\rm e}^{2i\pi \alpha (v-w)} \frac {\theta_1'(0)\theta_1(v-w -\mathcal X)}{\theta_1(\mathcal X)\theta_{1}(w-v)} \sqrt{\d v} \sqrt{\d w},\ \ \ \mathcal X:= \beta - \tau \alpha,
\ee 
where the multiplier system\footnote{We make here an abuse of notation where $\mathcal X$ denotes at the same time the point in the Jacobian and the corresponding character $\mathcal X :\pi_1(\mathcal C)\to \C^\times$.} is $\mathcal X(a) = {\rm e}^{2i\pi \alpha}$ and $\mathcal X(b) = {\rm e}^{2i\pi \beta}$, with $a,b$ denoting  the two cycles that correspond to $v\mapsto v+1$ and $v\mapsto v + \tau$, respectively. The Szeg\"o\ kernel consequently satisfies
\bea
S(v+1,w) ={\rm e}^{2i\pi \alpha} S(v,w),\ \ \
S(v,w+1) ={\rm e}^{-2i\pi \alpha} S(v,w),\nn\\
S(v+\tau,w) ={\rm e}^{2i\pi \beta} S(v,w),\ \ \
S(v,w+\tau) ={\rm e}^{-2i\pi \beta} S(v,w).
\eea
Following the general strategy, we can write a basis of 
$\scr P = \bigoplus_{n\geq 0} H^0\Big(\mathcal X\otimes \sqrt{\mathcal K} \big((n+1)\infty\big)\Big)$, $\scr P^\vee = \bigoplus_{n\geq 0} H^0\Big( \mathcal X^\vee\otimes \sqrt{\mathcal K} \big((n+1)\infty\big)\Big)$ as follows:
\bea
\varphi_\ell (w) =  \frac {\d^\ell }{\d w^\ell} \frac{S(v,w)}{\sqrt{\d w}}\bigg|_{w=0}
 \qquad
\varphi_\ell^\vee = \frac {(-1)^{\ell+1}\d^\ell }{\d w^\ell} \frac{S(w,v)}{\sqrt{\d w}}\bigg|_{w=0} , \ \ \ \ \ell=0,1,\dots
\eea
For example, the sections of $\scr P$ are given by linear combinations of the following form ($n\geq 1$)
\be
\varphi(v) = {\rm e}^{2i\pi v \alpha} \prod_{j=1}^n \frac {\theta_1(v-a_j)}{\theta_1(v)}\sqrt{\d v}
\ee
where   $a_1,\dots, a_n$ are any numbers
with $\sum a_j = \mathcal X = \beta - \alpha \tau \mod \Z+\tau \Z$. 
The Weyl sections are 
\bea
\mathcal W(v) &=  \sqrt{\d v}\int_{\gamma} {\rm e}^{2i\pi \alpha v} \frac {\theta_1'(0)^2\theta_1(v-w - \mathcal X)}{\theta_1(\mathcal X)^2\theta_{1}(w-v)} \frac{\theta_1\le(w- \mathcal X\ri)}
{\theta_1(w)}
{Y}(w)\d w\\
\mathcal W^\vee (v)& =  \sqrt{\d v}\int_{\gamma} {\rm e}^{-2i\pi \alpha v} \frac {\theta_1'(0)^2\theta_1(w-v - \mathcal X)}{\theta_1(\mathcal X)^2\theta_{1}(v-w)} \frac{\theta_1\le(w +  \mathcal X\ri)}
{\theta_1(w)}
{Y}(w)\d w
\eea
Consider the matrices 
\be
\label{defPhi}
\Phi(z):= \le[
\begin{array}{cc}
\varphi_0(v)  &\varphi_0(-v)\\
\varphi_1(v)   & \varphi_1(-v)
\end{array}
\ri], \ \ \Phi^\vee(z):= \le[
\begin{array}{cc}
\varphi_0^\vee(v)  & \varphi_1^\vee(v) \\
\varphi_0^\vee(-v)  & \varphi_1^\vee(-v)
\end{array}
\ri],\qquad z = \wp(v).
\ee
The following Lemma clarifies their mutual relationship:
\begin{lemma}
\label{lemmaphinv}
The inverse of $\Phi(z)$ is given by
\be
\Phi(z)^{-1} = \frac{-1} {4\omega_1^2 \d z}\Phi^\vee(z) \le[\begin{array}{cc}
0 & 1\\
1 & 0
\end{array}\ri]  = \frac{-1} {4\omega_1^2 \d z}  \le[
\begin{array}{cc}
 \varphi_1^\vee(v)  & \varphi_0^\vee(v) \\
\varphi_1^\vee(-v)  & \varphi_0^\vee(-v)
\end{array}
\ri].
\label{phinv}
\ee
\end{lemma}
{\bf Proof.}
The matrices $\Phi,\Phi^\vee$ satisfy 
\be
\Phi(z) \Phi^\vee(z) =-(2\omega_1)^2 \le[\begin{array}{cc}
0 & 1\\
1 & 0
\end{array}\ri] \d z.
\label{phiphivee}
\ee
This can be seen as follows; 
by a direct computation one has the expansions near $v=0$;
\bea
\varphi_0^{\{\vee\}}(v)& =   \frac 1{ v} (1 + \mathcal O(v))\sqrt{\d v}
, \qquad
\varphi_1^{\{\vee\}}(v) =  \frac 1{ v^2} (1 + \mathcal O(v^2))\sqrt{\d v}.\label{expphi1}
\eea
Then one needs to find the coefficient of $v^{-3}$ in the expansion of $\Phi\Phi^\vee$ at the origin. A direct computation yields 
\be
\Phi \Phi^\vee = 2 \le[
\begin{array}{cc}
0&1\\1&0\end{array}
\ri] \frac {\d v}{v^3}(1+ \mathcal O(v^2)).
\ee
 We then use that 
\be
\d v =   \frac {\d z}{2\omega_1 y}, \ \ \ y = (2\omega_1)^{-3}\wp'(v)\ ,\ \ \Rightarrow 
-2\omega_1^2\d z =  \frac {\d v}{v^3} (1+ \mathcal O(v)).
\ee 
So we conclude that \eqref{phiphivee} holds which implies that  the inverse of $\Phi$ is indeed given by \eqref{phinv}\QED

The Proposition \ref{propZ} states that for any section $\varphi\in \scr P_n$ we can write 
\be
\varphi(v) = A(z) \varphi_0(v) + B(z)\varphi_1(v),\ \ \ z= \wp(v)= Z.
\ee
with $A, B$ some polynomials in $z$.
It is seen directly that $\det \Phi$ is an elliptic function of $v$ with a triple pole at $v=0$ and zeros at $\frac 1 2, \frac \tau 2 , \frac {\tau+1}2$ and hence it is proportional to $\wp'(v)$. Then 
\be
A(z) = \frac 1{\det \Phi} \det \le[
\begin{array}{cc}
\varphi(v) & \varphi(-v) \\
\varphi_1(v)& \varphi_1(-v)
\end{array}
\ri], \ \ \ B(z) = \frac 1{\det \Phi} \det \le[
\begin{array}{cc}
\varphi_0(v)& \varphi_0(-v)\\
\varphi(v) & \varphi(-v) 
\end{array}
\ri].
\ee
Here $A,B$ are polynomials in $z=\wp(v)$ because the determinant in the numerator vanishes at the same half periods  (to order at least equal to) and the ratio is an even elliptic function with pole only at $v=0$. 
Unfortunately there are no known explicit orthogonal systems on elliptic curves. We consider some examples corresponding to finite matrix orthogonality that are relevant for applications originating from the study of the Aztec diamond \cite{Duits-Kuijlaars-aztec, KuijlaarsGroot} in the next section.

\subsection{$2R$--torsion points on elliptic curves and finite (matrix) orthogonality}
\label{sectorsion}
The following class of examples is prompted by the geometric interpretation of the nice example contained in \cite{Duits-Kuijlaars-aztec} studied in the context of the periodic tilings of the aztec diamond. We propose here a generalization of that setup. We comment in Remark \ref{compare} on how this is a generalization and how it compares to loc. cit.

Consider the  elliptic curve \eqref{ellyz}.
On its Jacobian $\mathbb J$, viewed as an abelian group, the $2R$--torsion points are the points $q\in \mathbb J$ such that $2R q\equiv 0$. Namely they are the points of the form
\be
q = \frac b {2R}+ \tau \frac a {2R},\ \ \ a,b\in \Z.
\label{tors}
\ee
Up to the equivalence $q\sim q+n+m\tau$ (for $m,n\in \Z$), we can assume $a,b\in\{0,\dots, 2R-1\}$. 
Another characterization of the $R$--torsion points which is relevant for us is the following: 
they are in correspondence with functions on the elliptic curve that have at most a  unique pole of order $R$,  a unique zero of the same multiplicity and such that both these points map to the same $z$--value (i.e. on the two sheets of the cover). 

Amongst those, the ones that are non-trivial (i.e. non-constant) correspond to the torsion points that are not periods or half-periods; moreover any non-trivial such function ${Y}$ has the property that the elliptic involution maps ${Y}\to {Y}^{-1}$ and hence we want to consider these torsion points modulo the elliptic involution as well.

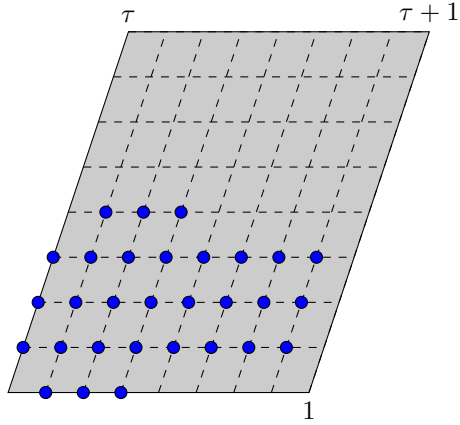
\begin{figure}
\begin{center}
\begin{tikzpicture}[scale=4]
\pgfmathsetmacro{\R}{4} ;
\pgfmathsetmacro{\RR}{int(\R+\R)} ;
\pgfmathsetmacro{\Rd}{int(\R-1)} ;
\coordinate (tau) at (0.4,1.2);
\draw [fill=black!20!white] (0,0) to (tau)  to ($(tau)+ (1,0)$) to (1,0) to cycle;
\foreach  \x in {1,...,\RR} {
\draw [dashed](\x/\RR ,0) to +(tau);
\draw [dashed]($\x/\RR*(tau)$) to +(1,0);
};
\foreach \x in {1,...,\Rd}{
\draw[fill=blue] ($\x/\RR*(tau)$) circle [radius=0.02];
\draw[fill=blue] ($1/2*(1,0)+ \x/\RR*(tau)$) circle [radius=0.02];
\draw[fill=blue] ($\x/\RR*(1,0)$) circle [radius=0.02];
\draw[fill=blue] ($\x/\RR*(1,0) + 1/2*(tau)$) circle [radius=0.02];
};
\foreach \x in {1,...,\Rd}{
\foreach \y in {1,...,\Rd}{
\draw[fill=blue] ($\x/\RR*(tau) + \y/\RR*(1,0)$) circle [radius=0.02];
\draw[fill=blue] ($(1/2,0)+\x/\RR*(tau) + \y/\RR*(1,0)$) circle [radius=0.02];
};
};
\node [above]at (tau) {$\tau$};
\node [above]at ($(tau)+(1,0)$) {$\tau+1$};
\node [below]at (1,0) {$1$};
\end{tikzpicture}
\end{center}
\caption{The $2R$-torsion points (here depicted for $R= 4$) are the points at the intersection of the dashed grid, modulo $\Z+ \tau \Z$. The points that correspond to our functions are the torsion points that are not half-periods and modulo the hyperelliptic involution $v\to -v$. They are  indicated by dots. Their total number is $2R^2-2$. 
}
\label{figtorsion}
\end{figure}

The $2R$--torsion points are actually algebraic points on the curve (meaning that  they can be written as the solution of an algebraic equation in the coefficients of the curve). 
To find them one writes the most general elliptic function with pole of order $R$ at $z_\star$:
\be
F:= \frac {p_R(z) - y(z) q_{R-2}(z)}{(z-z_\star)^R}
\ee
where $p,q$ are polynomials of the indicated degrees. 
The condition is that the numerator, $N(z)$, satisfies $N^{(\ell)}(z_\star)=0$ for $\ell = 0,1,\dots, 2R-1$ with the determination of $y(z)$ at $z_\star$ being the one corresponding to the chosen sheet. For $F$ to be non-constant we need to have $q_{R-2} \not\equiv0$. The equations $N^{(R+\ell)}(z_\star)=0$ for $\ell=1,\dots, R-1$ give a linear homogeneous system for the coefficients of $q_{R-1}$ and hence the condition is that the determinant of this system vanishes. 
A direct computation shows that 
\be
\det \le[
\begin{array}{cccccc}
y^{(R+1)} & (R+1) y^{(R)} & R(R+1) y^{(R-1)} & \dots  & \frac {(R+1)!}{(R-3)!} y^{(3)}\\
y^{(R+2)} & (R+2) y^{(R+1)} & (R+1)(R+2) y^{(R)} & \dots  & \frac {(R+1)!}{(R-2)!} y^{(4)}\\
\vdots &&\vdots\\
y^{(2R-1)} & (2R-1) y^{(R)} & (2R-2)(2R-1) y^{(2R-1)} & \dots  & \frac {(2R-1)!}{(2R-3)!} y^{(R+1)}
\end{array}
\ri]=0
\ee
where all derivatives are evaluated at $z_\star$. This gives a polynomial equation of degree $2 R^2-2$ in $z_\star$ and this is the number of (inequivalent) $2R$ torsion points. See Fig. \ref{figtorsion}.

Here is how these functions can be expressed using the Jacobi theta function (up to multiplicative scalar):
\be
\label{{Y}R}
{Y}(v) = \le(
{\rm e}^{2i\pi \frac aR v} \frac {\theta_1\le( v + \frac{b+a\tau}{2R}\ri)}{\theta_1\le(v-\frac{b+a\tau}{2R}\ri)}
\ri)^R
\ee
They also have the evident property that ${Y}(-v)=\frac 1{{Y}(v)}$, namely, the elliptic involution that exchanges the sheets of the $z$--projection acts as reciprocal on these functions. This latter condition fixes the multiplicative scalar to be $\pm 1$. 

\begin{remark}[``Prime'' torsion points]
\label{remprime}
Note that the expression in the bracket in \eqref{{Y}R} is {\it not} an elliptic function but rather a meromorphic  section of a flat bundle $\mathcal X_{a,b}$. This is the bundle with multiplier system 
${\rm e}^{\frac{2i\pi a}{R}}$, ${\rm e}^{-\frac{2i\pi b}{R}}$ on the $A$, $B$ cycles, respectively.
Namely,  although it has pole and zero of multiplicity $R$, the function ${Y}(v)$ is not the $R$ power of a single--valued function with simple pole and zero, which cannot possibly exist on a curve of genus $g\geq 1$. 

However, for example, if $R$ is not a prime number and $a,b$ have the same common divisor with $R$, then some of these functions are actually powers of that common divisor. 
We call ``prime'' those  torsion points \eqref{tors} where at least one of $a,b$ is prime relative to $2R$. 
\end{remark}
\begin{remark}
The case $R=1$ of formula \eqref{{Y}R} gives constant functions. The case $R=2$ is the first interesting case, corresponding to quarter periods. This is the case implicitly used in \cite{Duits-Kuijlaars-aztec} (see Sec. \ref{compare}).
\end{remark}

Let us set $z_\star:= \wp(\frac {b+a\tau}{2R})$. According to our notation $z_\star^{(j)}$, $j=1,2$ are the two points on the elliptic curve corresponding to $\wp(v) = z_\star$ and hence they are (in the $v$--plane) $\pm \frac { b+a\tau}{2R}$. We now define 
\bea
\psi_{2R} (v):=\frac{(z-z_\star)^R}{{Y}(v)}\varphi_0(v),\qquad
\psi_{2R+1} (v):=\frac{(z-z_\star)^R}{{Y}(v)}\varphi_1(v),\ \ \ z=\wp(v),
\eea
with similar expressions for $\psi_{2R}^\vee, \psi_{2R+1}^\vee$. Note that all these sections have a zero of order at least $2R$ at $z_\star^{(1)}$ on one the main sheet, and are analytic at $z_\star^{(2)}$.

We define now the pairing \eqref{pairing} in this case to be 
\be
\label{pairingg1}
\langle\varphi, \varphi^\vee\rangle = \oint_{|v-z_\star^{(1)}|=\epsilon} \varphi(v) \varphi^\vee(v) {Y}^2(v),\qquad \varphi\in \scr P,\ \ \varphi^\vee\in \scr P^\vee.
\ee
Note that we may integrate also on a circle around $z_\star^{(2)}$ but it would yield zero because there is no pole of the integrand there.

It is evident, since ${Y}$ has a pole of order $R$, that this pairing is of rank $2R$ and $\psi_{2R, 2R+1}$ (and their $\vee$ counterparts) span  the $2$--dimensional kernel of the pairing in $\scr P_{2R+1}, \scr P_{2R+1}^\vee$, respectively.

\paragraph{Matrix biorthogonality.}
We define, according to the general picture, the weight matrix $W(z)$  by the formula
\be
\sqrt{W(z) } = \Phi(z) \le[
\begin{array}{cc}
{Y}(z^{(1)}) & 0\\ 0 & {Y}(z^{(2)})
\end{array}\ri]\Phi(z)^{-1},
\label{defW}
\ee
with $\Phi$ as in \eqref{defPhi} and $\Phi^{-1}$ given as in Lemma \ref{lemmaphinv}. 
Note that $\sqrt{W}$  is a rational function of $z$ with a pole of order $R$ at $z=z_\star$ and with unit determinant.
Moreover it is bounded as $z\to \infty$. 
The matrix
\be
\label{PsiR}
\Psi_R(z) = \le[
\begin{array}{cc}
\psi_{2R}(v)  & \psi_{2R}(-v)  \\
\psi_{2R+1}(v)  & \psi_{2R+1}(-v)  
\end{array}
\ri],\ \ \ z = \wp (v),
\ee
can then be written as per Prop. \ref{propZ}
\bea
\Psi_R(z) &= P_R(z) \Phi(z),\\
P_R(z) &= (z-z_\star)^R \Psi_{2R}(z) \Phi^{-1}(z) = 
(z-z_K)^{R} \Phi(z) \le[
\begin{array}{cc}
{Y}^{-1}(v) & 0 \\
0 & {Y}^{-1} (-v)
\end{array}
\ri] \Phi^{-1}(z)=\cr
&= (z-z_\star)^R \sqrt{W^{-1}}(z)
\eea
By construction, or by inspection, one concludes with the following 
\begin{lemma}
The polynomial matrix $P_R$ of degree $R$ satisfies 
\be
\oint_{|z-z_\star|=\epsilon} P_R(z) W(z) Q(z)\d z =0
\ee
for all polynomial matrices $Q(z)$. 
\end{lemma}

Following the ideas in \cite{Duits-Kuijlaars-aztec, KuijlaarsGroot} can find explicitly the matrix polynomial $P_{R-1}$ which is biorthogonal to all powers of $z^j\1$, $j=0,\dots, R-2$. 

Interestingly, following the idea in \cite{Duits-Kuijlaars-aztec}, one can write the MOP of degree $R-1$ in a rather explicit form as in the  following Lemma.
\begin{lemma}
\label{lemmaPR-1}
The matrix--valued function  
\be
\label{PR-1}
P_{R-1}(z) = (z-z_\star)^{R} \oint_{|w-z_\star|=\epsilon}  \Psi_0(w) \mathbb E_{11}\Psi_0^{-1}(w)\frac{ \d w}{(w-z_\star)^{2R}(w-z)2i\pi} \sqrt{W^{-1}}(z)
\ee
is a polynomial of degree $R-1$. 
\end{lemma}
{\bf Proof.}
Since $\sqrt{W}$ is a rational matrix and bounded at infinity, it is clear from the expression that $P_{R-1}$ is of growth bounded by $z^{R-1}$ at infinity. The residue at $z=z_\star$ is a  Laurent polynomial in $(z-z_\star)^{-1}$ of degree at most $2R-1$. We now show that the total expression \eqref{PR-1} does not have a pole at $z_\star$. To see this we use Cauchy's residue theorem on the integral,  taking a circle centered at $z_\star$ that contains $z$. Then we find 
\bea
P_{R-1}(z) = (z-z_\star)^{R} \oint_{|w-z_\star|>|z-z_\star|}  \Psi_0(w) \mathbb E_{11}\Psi_0^{-1}(w)\frac{ \d w}{(w-z_\star)^{2R}(w-z)2i\pi} \sqrt{W^{-1}}(z)
\nn
\\
-  \Psi_0(z) \mathbb E_{11}\Psi_0^{-1}(z)\frac{ 1}{(z-z_\star)^{R}} \sqrt{W^{-1}}(z).
\label{1312}
\eea
The term on the first line of \eqref{1312} is analytic at $z=z_\star$ because $\sqrt{W^{-1}}(z)(z-z_\star)^R$ is analytic and so is the integral. 
The term on the second line of \eqref{1312} may still, in principle, have a pole of order $R$.  However, by the definition \eqref{defW} we have that it matches 
\be
 \Psi_0(z) \mathbb E_{11}\Psi_0^{-1}(z)\frac{ 1}{(z-z_\star)^{R}} \sqrt{W^{-1}}(z)=
\Psi_0(z) \le[
\begin{array}{cc}
\frac{1}{(z-z_\star)^R {Y}(z^{(1)} )}& 0\\ 0 &0
\end{array}\ri]\Psi_0(z)^{-1}.
\ee 
Since ${Y}$ has a pole at $z^{(1)}$ of order $R$, we see that $\frac{1}{(z-z_\star)^R {Y}(z^{(1)} )}$ is analytic at $z=z_\star$. 
\QED

We now show that $P_{R-1}$ is also a biorthogonal matrix polynomial.
\bp
\label{propPR-1}
The polynomial $P_{R-1}$ is orthogonal to all matrix polynomials of degree $\leq R-2$:
\be
\oint_{|z-z_\star|=\epsilon} \hspace{-15pt} \d z\,\,\, P_{R-1}(z) W(z) z^\ell =0,\ \  \ \ \ell=0,\dots, R-2.
\ee
\ep
{\bf Proof.}
We compute directly, 
\bea
\oint_{|z-z_\star|=2\epsilon} \hspace{-14pt} \d z \,\,P_{R-1}(z) W(z) z^{j} =
\oint_{|z-z_\star|=2\epsilon}\hspace{-13pt}\d z  \oint_{|w-z_\star|=\epsilon}\hspace{-13pt} \d w\,\,  \Psi_0(w) \mathbb E_{11}\Psi_0^{-1}(w)\frac{ z^j
(z-z_\star)^{R}  \sqrt{W}(z)
}{(w-z_\star)^{2R}(w-z)2i\pi}
\label{320}
\eea
By Fubini's theorem we can exchange the order of integration. Then, using the Cauchy theorem we can  replace the integration on $|z-z_\star|=2\epsilon$ by an integration on $|z-z_\star=\frac \epsilon 2$ while also picking up the residue at $z=w$:
\bea
\eqref{320} = 
  \oint_{|w-z_\star|=\epsilon}\hspace{-13pt} \d w\,\,\oint_{|z-z_\star|=\frac \epsilon2} \hspace{-14pt} \d z   
  \Psi_0(w) \mathbb E_{11}\Psi_0^{-1}(w)\frac{ z^j
(z-z_\star)^{R}  \sqrt{W}(z)
}{(w-z_\star)^{2R}(w-z)2i\pi}
\nn\\
- \res{w=z_\star} \Psi_0(w) \mathbb E_{11}\Psi_0^{-1}(w)\frac{ w^j
  \sqrt{W}(w)\d w
}{(w-z_\star)^{R}}
\label{319}.
\eea
The term on the first line of \eqref{319} vanishes because $(z-z_\star)^R\sqrt W$ is locally analytic. We now show that also the residue on the second line vanishes.  Indeed we are computing the residue of 
\be
w^j\Psi_0(w) \le[
\begin{array}{cc}
\frac{ {Y}(w^{(1)} )}{(w-z_\star)^R}& 0\\ 0 &0
\end{array}\ri]\Psi_0(w)^{-1}.
\ee
The residue is the same as the residue of 
\be
w^j\Psi_0(w) \le[
\begin{array}{cc}
\frac{ {Y}(w^{(1)} )}{(w-z_\star)^R}& 0\\ 0 &\frac{ {Y}(w^{(2)} )}{(w-z_\star)^R}
\end{array}\ri]\Psi_0(w)^{-1} = \frac {w^j}{(w-z_\star)^R} \sqrt W(w),\label{321}
\ee
because $\frac{ {Y}(w^{(2)} )}{(w-z_\star)^R}$ is locally analytic at $z_\star$ thanks to the fact that ${Y}$ has a {\it zero} of order $R$ at $z_\star^{(2)}$.  But now, $\sqrt{W}(w)$ is a {\it rational} matrix, bounded at infinity. So, the residue of the expression \eqref{321} can be computed as the residue at infinity. The latter is zero if $j\leq R-2$. 
\QED
\begin{example}
For the family of curves $y^2 = z(z-a)(z-a^{-1})$ the points $z = \frac{\pm 1  + \sqrt{1-a^2}}a, a\pm \sqrt{a^2-1}$ are quarter periods.

For $y^2 = z(z-1)(z-8+4\sqrt{3})$ one of the $16$ sixth periods ($R=3$) is $z_\star = 2$. 
\end{example}
\subsection{A more general example in genus 1.}
\label{moregeneral}

Riding the same train of thoughts, we can consider for weight the elliptic functions having $R$ distinct pairs of simple zeros/poles at points $z_{j}$ instead of \eqref{{Y}R}:
\be
{Y}(v):= {\rm e}^{2i\pi a v}\prod_{j=1}^R \frac {\theta_1( v + v_j)}{\theta_1(v-v_j)} 
\ee
where 
\be
\sum_{j=1}^R v_j = \frac{b+\tau a}2, \ \ a,b\in \Z. \label{Rhalf}
\ee
The function has simple poles at $z_j:= \wp (v_j),\ j=1,\dots, R$ on one sheet and simple zeros at the same points on the other sheet; the values of these $z_j$ are not independent because of the constraint \eqref{Rhalf}. 
\begin{remark}
To clarify the constraint alluded at above, we note that the function ${Y}$ can be written as 
\be
{Y}(p) = \frac{ A(z) - B(z) y(z)}{\prod_{j=1}^R (z-z_j)},
\ee
where $A, B$ are polynomials of degrees $R, R-2$, respectively, with $B$ not identically zero and $A$ monic (so as to dispose of the freedom of multiplication by a scalar). 
Then the constraints are that the numerator and its derivative must vanish at $z_j$:
\be
A(z_j) - B(z_j) y(z_j)=0,\  \  \ 
A'(z_j) - B'(z_j) y(z_j) - B(z_j) y'(z_j) =0.
\ee
This gives $2R$ equations for the coefficients of $A$ (monic), $B$ and the $R$ positions of the points $z_j$. 
So we see that the solutions form a family of dimension $(R+R-1+R)-2R = R-1$.  \hfill $\triangle$
\end{remark}
Define 
\bea
\psi_{2R} (v):=\frac{q_r(z)}{{Y}(v)}\varphi_0(v),\qquad
\psi_{2R+1} (v):=\frac{q_r(z)}{{Y}(v)}\varphi_1(v),\ \ \ z=\wp(v),\  \ \ \ q_R:= \prod_{j=1}^R(z-z_j),
\eea
with similar expressions for $\psi_{2R}^\vee, \psi_{2R+1}^\vee$. Once again, $\psi_{2R,2R+1}$ (and the $\vee$ versions) have zeros at $z_j$'s on one sheet and are holomorphic at the same points on the second sheet.  
We then define the weight matrix $W$ by the same formula  \eqref{defW}.
Observe that $\sqrt{W}$ is still a rational function of $z$ with unit determinant and with poles at the $z_j$'s and bounded at $z=\infty$ ($v=0$).
The pairing is defined by the same \eqref{pairingg1} with the contour of integration taken as the union of small circles around $z=z_j$'s (or a loop containing all of them at once). In the $v$--plane, we have  $R$-pairs of circles surrounding $v=\pm v_j$ (modulo $\Z+\tau \Z$) but the circles around $-v_j$ do not contribute to the pairing. 
The matrix $\Psi_R$ is defined by the same \eqref{PsiR} but now the matrix-polynomial $P_R$ is of the form:
\bea
P_R(z)&= q_R(z) \sqrt{W^{-1}}(z)
\eea
In lieu of Lemma \ref{lemmaPR-1} and formula \eqref{PR-1} we have the formula 
\be
P_{R-1}(z) = q_R \oint_{|w-z_\star|=\epsilon}  \Psi_0(w) \mathbb E_{11}\Psi_0^{-1}(w)\frac{ \d w}{q_R(w)^2(w-z)2i\pi} \sqrt{W^{-1}}(z).
\ee 
The proof that $P_{R-1}$ is indeed a polynomial proceeds in the same way as for Lemma \ref{lemmaPR-1}. The orthogonality to all powers of $z^j$, $j=0,\dots, R-2$ is also shown as in Prop. \ref{propPR-1}, simply replacing expressions of the form $(z-z_\star)^R$ by $q_R(z)$.

\begin{remark}[Comparison with \cite{Duits-Kuijlaars-aztec}]
 \label{compare}
In loc. cit. the authors consider the weight matrix (in their paper $N=2k$)
\be
W(z) = W_1(z)^{2k} ,\ \ \ 
W_1(z):= \frac {1}{(z-1)^2}\le[
\begin{array}{cc}
(z+1)^2 + 4\a^2 z & 2\a(\a + \a^{-1}) (z+1)\\
2\a^{-1}(\a+\a^{-1}) z(z+1) & (z+1)^2 + 4\a^{-2} z
\end{array}
\ri], \ \  \a >0.
\ee
Note that $\det W_1=1$ and $\sqrt{W} = W_1^k$ is a rational matrix as in our description. 
The spectral analysis of $W_1$  is as follows. The spectral curve $\mathcal C$ is of genus $1$ and given by  (see (2.11) in loc. cit.) 
\be
\label{DKcurve}
y^2 = z(z+\a^2)(z+\a^{-2}).
\ee
The quarter-period points of the curve \eqref{DKcurve} are $z_\star = \pm 1, \a^2\pm \sqrt{\a^4-1}, 
-\a^{-2} \pm \sqrt{\a^{-4}-1}$. In particular, note that the pole is at the quarter period $z_\star =1$. 

The eigenvalue of $W_1$ on the spectral curve is the function 
$$
{Y}_1 = \frac { 2(a^4+1)z +a^2(z+1)^2 - 2a(a^2+1)y }{a^2(z-1)^2}
$$
which can be seen to have a double zero  at $(z_\star,y_\star) = (1, \sqrt{(1+\a^2)(1+\a^{-2})}) $ and double pole at $(z_\star, -y_\star)$. This is a reflection of the fact that $z_\star =1$ is a quarter period as explained in the previous section.
The weight $W(z)$ has eigenvalue ${Y}^2 = {Y}_1^{2k}$ which has pole of order $4k$ and hence $R=2k$ in the notation of Section \ref{sectorsion}.  We see here that  the function ${Y}$  is a perfect $k$ power , i.e. does not correspond to ``prime'' torsion points (in the terminology of Rem. \ref{remprime}).

The divisor for the eigenvector line bundle $\scr D$ of  Sec. \ref{abelianization} is obtained by looking at the poles of the eigenvector normalized to have $1$ in the first entry and then found to be $\scr D = (0) + (-1)^{(1)}$. This divisor is equivalent to $(\infty) + (-1)^{(2)}$ which is seen by considering the function 
$$
f(z,y):= \frac {z c -y}{z(z+1)} ,\ \ \qquad c = \sqrt{(\a^2-1)(1-\a^{-2})}= \frac{\a^2-1}{\a} 
$$
which has divisor of poles at $(0)+(-1)^{(1)}$ and divisor of zeros at $(\infty) + (-1)^{(2)}$. 
So the multiplier system $\mathcal X$ corresponds to the quarter period $(-1)^{(2)}$. 
\end{remark}

\paragraph*{Acknowledgements.}
The work  was supported in part by the Natural Sciences and Engineering Research Council of Canada (NSERC) grant RGPIN-2016-06660.

\appendix
\renewcommand{\theequation}{\Alph{section}.\arabic{equation}}

\section{Theta functions and the Szeg\"o\ kernel}
\label{theta}
We choose a Torelli marking $\{a_1,\dots, a_g, b_1,\dots, b_g\}$ for $\mathcal C$ in terms of which we construct the classical Riemann Theta functions. We refer to \cite{Fay}, Ch. 1-2 for a review of these classical notions.  The corresponding normalized Abelian differentials, period matrix and Abel map will be denoted by $\omega_j, \bs \tau, \mathfrak A$, respectively:
\be
\oint_{a_j} \omega_k = \delta _{jk},\  \ \bs \tau_{ij}= \oint_{b_j} \omega_i = \oint_{b_i} \omega_j,\ \ \ \mathfrak A_{p_0}(p)= \int_{p_0}^p \vec \omega: \mathcal C\times \mathcal C \to \mathbb J(\mathcal C).
\ee
In our setting the basepoint, $p_0$, of the Abel map will be chosen amongst the support of $\Xi= Z^{-1} (\infty)$. The basepoint will not be indicated except when necessary and the Abel map is naturally extended to a map on divisors. 

For given characteristics $[\bs\alpha, \bs \beta]\in \R^g\times \R^g$ corresponding to the point  ${\mathbf e}=\bs \beta - \bs \tau \bs \alpha\in \mathbb J(\mathcal C)$ in the Jacobian we denote with $S(p,q)= S_{[\bs\alpha, \bs \beta]}(p,q)$ the {\it Szeg\"o} kernel:
\be
S(p,q)= S_{[\bs\alpha, \bs \beta]}(p,q) =\frac{ \Theta_{[\bs\alpha, \bs \beta]}(p-q)}{\Theta_{[\bs\alpha, \bs \beta]}(0) E(p,q)}\label{Szegodef}
\ee
where $E(p,q)$ is the Klein prime form [\cite{Fay}, pag. 16] and $\Theta_{[\bs\alpha, \bs \beta]}$ denotes the Riemann Theta function with characteristics (see loc. cit.). 
Following tradition, the Abel map of points or divisors used in  arguments of Theta functions is understood for brevity.

The celebrated Fay identity, generalizing the tri-secant one, is the following equality 
\be
\label{Fayid}
\frac {
\Theta_{[\bs\alpha, \bs \beta]} \le(\sum_{j=1}^K (p_j- q_j)\ri)
} {\Theta_{[\bs\alpha, \bs \beta]}(0)}
\frac{\prod_{j<k} E(p_j,p_k) E(q_j, q_k)}  {\prod_{j,k} E(p_j, q_k)} = \det \bigg[S(p_a,q_b)\bigg]_{a,b=1}^K
\ee
which is valid as long as $\mathbf e:= \bs \beta- \bs \tau \bs \alpha$ does not belong to the Theta divisor: $\Theta(\mathbf e)\neq 0$. Points ${\mathbf e}\in \mathbb J(\mathcal C)$ are in bijective correspondence with unitary characters $\mathcal X:\pi_1\to U(1)$ as follows: if $a_j, b_j$ form a canonical symplectic basis in homotopy (which we can think also as a symplectic basis in homology), then one sets
\be
\mathcal X(a_j)= {\rm e}^{2i\pi \alpha_j},\ \ \ \ \mathcal X(b_j)= {\rm e}^{2i\pi \beta_j},\ \ j=1,\dots, g,
\ee
where $\alpha_j, \beta_j$ are the (real) characteristics of the point ${\mathbf e}$. Then the condition $\Theta(\mathbf e)=0$ has the geometric interpretation that  $h^0(\mathcal X \otimes \sqrt{\mathcal K})\geq 1$. Furthermore the order of vanishing of $\Theta$ at $\mathbf e$ is precisely the value of $h^0(\mathcal X \otimes \sqrt{\mathcal K})$.   

Let $\scr Q:= \sum_{j=1}^K q_j$. Then  the identity \eqref{Fayid} can be interpreted in the same way as the Vandermonde determinant; namely, it states whether the $K$ spanning sections $\varphi_\ell(p) := S(p,q_\ell)$ of  $H^0\big(\mathcal X\otimes \sqrt{\mathcal K}(\scr Q)\big)$ can interpolate any local section defined  at the points of the divisor $\scr P  = \sum_j p_j$. 

A small technical issue is when the divisor $\scr Q$ has points of higher multiplicity. As an extreme example if $\scr Q= K q_1$ then one proceeds as follows to obtain holomorphic sections of $\mathcal X \otimes \sqrt{\mathcal K}(\scr Q)$;
fixing a local coordinate $\zeta$ near $q_1$ ($\zeta(q_1)=0$)  one can write, for $q$ is the coordinate patch
\be
S(p,q) = \varphi(p,\zeta) \sqrt \d \zeta
\ee
where $\varphi(p,\zeta)$ is a section of the bundle $\mathcal X \otimes \sqrt K$ with a simple pole at $p=q$ (or, equivalently, a {\it holomorphic} section of $\mathcal X \otimes \sqrt K(q)$). Then one obtains sections with higher order poles at $q_1$ simply by 
\be
\varphi_j(p):= \frac {\d^{j-1}}{\d \zeta^{j-1}} \varphi(p,\zeta) \bigg|_{\zeta=0}, \ \ j=1,\dots, K
\ee
which can be interpreted as {\it holomorphic} spanning sections of $\mathcal X\otimes \sqrt{\mathcal K}(Kq_1)$
In this case the Fay identity \ref{Fayid} survives by taking appropriate limits and using l'Hopital's rule and gives
\be
\label{Fayiddeg}
\frac {
\Theta_{[\bs\alpha, \bs \beta]} \le(\sum_{j=1}^K p_j - Kq_1\ri)
} {\Theta_{[\bs\alpha, \bs \beta]}(0)}
\frac{\prod_{j<k} E(p_j,p_k)}  {\prod_{j} E(p_j, q_1)^K \d \zeta(q_1)^\frac K2} = \det \bigg[\varphi_a(p_b)\bigg]_{a,b=1}^K.
\ee
Here the expression $E(p,q_1) \sqrt {\d \zeta(q_1)}$ simply means the prime form where we have evaluated it  using the local coordinate $\zeta$ with respect to the variable $q_1$. 
\subsection{Proof of Theorem \ref{Propideal}}
\label{proof}
We  take $K= r=\deg Z$ and specialize the  formula  \eqref{Fayid} with $p_j = z^{(j)}$ and $\scr Q= \Xi$.
We choose the basis $\varphi_j$ described before Proposition \ref{propPhiPhivee},  and expressed in terms of the Szeg\"o\ kernel (or appropriate derivatives in the case $\Xi$ is not a simple divisor). 
 Then the argument of $\Theta_{[\bs\alpha, \bs \beta]}$ in the numerator of the left side of the Fay identity is zero because the evaluation of Abelian differentials at the points of a rational map yields the pullback of a holomorphic differential on $\C \P^1$ which is necessarily zero.
 
This is due to the fact that if $\omega$ is a holomorphic differential on $\mathcal C$ and we evaluate it at the $r$ preimages of $z$ then 
$
\sum_{j=1}^r \omega(z^{(j)})=0, \ \ \ \forall z\in \C.
$
Applying this to the basis of holomorphic differential and integrating, we deduce that 
$$
\sum_{j=1}^r \mathfrak A(z^{(j)}) 
$$
is a constant (independent of $z$). Thus 
\be
\Theta_{[\bs\alpha, \bs \beta]} \le(\sum_{j=1}^r; (z^{(j)} - \infty^{(j)}) \ri) = \Theta_{[\bs\alpha, \bs \beta]} \le(0\ri).
\ee 
Thus the Fay identity \eqref{Fayid} (or the likes of it as \eqref{Fayiddeg} if $\Xi$ has multiplicities) shows that the determinant, as a function of $z$,  vanishes only at the critical points of $Z$. 

Let $c$ be a critical value of $Z$ and $p_\ell$ the points above it. If $p_\ell$ has multiplicity $\nu_\ell$ then a local coordinate is $\zeta_\ell=(z-c)^{\frac 1 {\nu_\ell}}$. 

We trivialize the bundle $\mathcal X \otimes \sqrt {\mathcal K}$ near the points $Z^{-1}(c)$ using the above coordinates; this simply means that we write every section $\varphi_j(p)= f_j(\zeta_\ell) \sqrt{\d \zeta_\ell}$ and similarly $\psi_j(p) = g_j(\zeta_\ell) \sqrt{ \d \zeta_\ell}$ where $f_j, g_j$ are analytic functions near $0$ in the respective variables. 

The function $F$ is thus the ratio of the determinants of two  alternant matrices of the sets of functions $f_j, g_j$.  

 Near $z=c$ the term responsible for the vanishing of $\det \Phi$ is the Vandermonde--like term $\prod_{j<k} E(z^{(j)}, z^{(k)})$ in the formula \ref{Fayid} (or the likes of \eqref{Fayiddeg} in case $\Xi$ is not a simple divisor). Letting $\mathbb E(\zeta, \xi) = E(p,q) \sqrt{ \d \zeta(p) \d \xi(q)}$ be the local trivialization of the Klein prime form in the coordinates $\xi,\zeta$ we have 
\be
\prod_{j<k} \mathbb E(z^{(j)}, z^{(k)}) \propto \prod_{\ell=1}^{|Z^{-1}(c)|} (z-c)^{\frac {\nu_\ell-1}{\nu_\ell}}.
\ee
Now we show that the determinant of  the alternant matrix in the numerator of $F$ in \eqref{F} must vanish at least of the same order. 

To simplify the matter and without major loss of generality, we assume that there is only one branchpoint, $p^*$, above $z=c$ of multiplicity $\nu$, with $\zeta = (z-c)^{\frac 1{\nu}}$ being a local coordinate. We also enumerate the sheets of the map $Z$ so that the first $\nu$ sheets are glued at $z=c$.
The functions $g_j$ near $p^*$ are Puiseux series
\be
g_j(p) = g_j^{0}  +  g_j^1 (Z(p)-c)^\frac 1 \nu + \dots = g_j^0 + g_j^1 \zeta + \mathcal O(\zeta^2).
\ee
The different sheets are distinguished by multiplying $\zeta$ by a $\nu$--root of unity;
\be
g_j(z^{(a)}) = g_j^0 + g_j^1 {\rm e}^{\frac{2i\pi a}{\nu}} \zeta + \mathcal O(\zeta^2), \ \ \ \ a=1,\dots, \nu.
\ee
Thus the alternant matrix has the following expansion 
\be
\le[
\begin{array}{c|c c}
G^0 & F^0
\end{array}
\ri] + \zeta \le[
\begin{array}{c|c}
\Omega  & 0_{r\times (r-\nu)}
\end{array}
\ri]  + \mathcal O(\zeta^2)
\ee
where $\Omega\in \mathrm{Mat}_{r\times \nu}$ is the matrix $\Omega_{ja}= g_j^{1} {\rm e}^{\frac {2i\pi a}{\nu}}$ and $G^0\in \mathrm{Mat}_{r\times \nu}$ is the matrix $G^0_{ja} = g_j^0$  of rank $1$ (all columns are equal). The matrix $F\in \mathrm{Mat}_{r\times (r-\nu)}$ is in general of maximal rank $r-\nu$ and it is irrelevant for the present consideration (it consists of the evaluations of the functions $g_j$ at the distinct points $ z^{(a)}, \ a=\nu+1,\dots, r$).

It is then clear that in the evaluation of the determinant the lowest order term are of order $\zeta^{\nu-1}= (z-c)^{\frac {\nu-1}\nu}$. 
This proves that the ratio $F(z)$ is bounded at $z=c$; since it is single--valued, the singularity is removable and $F(z)$ extends to an analytic function at $z=c$. The case where there are several branchpoints in the fiber of the critical value $z=c$ is only slightly more complicated and left to the reader.
\QED
\begin{remark}
The arguments above apply also if  $Z$  maps to a higher genus curve; in that case the ratio $F$ will be simply a meromorphic function on the target curve.
\end{remark}

\end{document}